\documentclass[a4paper,10pt,twoside,leqno]{article}
\pdfoutput=1
\usepackage[utf8]{inputenc}
\usepackage[T1]{fontenc}
\usepackage[full]{textcomp}

\usepackage{csquotes}
\usepackage[english]{babel}

\usepackage{mathtools}
\mathtoolsset{mathic} 

\usepackage{cfr-lm}
\usepackage{dsfont}  

\usepackage{microtype}
\usepackage{ifpdf}
\ifpdf
\def\letterspace{\textls}
\else
\def\letterspace{}
\fi
\usepackage[
nohead,nomarginpar,
paper=a4paper,
]{geometry}

\usepackage{lastpage}
\usepackage{fancyhdr}
\usepackage{longtable}
\usepackage{booktabs}
\usepackage{tabu}

\usepackage[inline]{enumitem}
\setlist{noitemsep,nosep,listparindent=\parindent}
\setlist[itemize]{label=\guillemotright}
\setlist[enumerate,1]{ref=\thesubsection.\arabic*}
\setlist[enumerate,2]{label=\alph*.,ref=\theenumi.\alph*}

\usepackage[backend=biber,doi=false,url=false,isbn=false,%
safeinputenc,style=quasialphabetic,citestyle=alphabetic,%
]{biblatex}
\bibliography{\jobname.bib}
\defbibheading{bibliography}[\bibname]{\sectionstar{#1}}
\usepackage{ifthen}
\makeatletter
\renewcommand{\part}[1]{%
 \cleardoublepage%
 \vbox{\null\vskip90pt%
 \normalfont\fontsize{20pt}{30pt}\selectfont%
 \baselineskip=30pt%
 \scshape\noindent\letterspace{#1}\par}%
 \addcontentsline{toc}{part}{#1}%
 \@afterindentfalse%
 \@afterheading%
}
\renewcommand{\section}[1]{%
 \refstepcounter{section}%
 \vbox{\null\vskip15pt%
 \normalfont\fontsize{12pt}{15pt}\selectfont%
 \baselineskip=15pt%
 \centering\scshape\noindent\letterspace{\thesection\quad#1}%
 \par}
 \addcontentsline{toc}{section}{\protect\numberline{\thesection} #1}%
 \@afterindentfalse%
 \@afterheading%
} 
\newcommand{\sectionstar}[1]{%
 \vbox{\null\vskip15pt%
 \normalfont\fontsize{12pt}{15pt}\selectfont%
 \baselineskip=15pt%
 \centering\scshape\noindent\letterspace{#1}%
 \par\nobreak\vskip15pt}
 \@afterindentfalse%
 \@afterheading%
} 
\renewcommand{\paragraph}[1]{\par\bigskip\refstepcounter{subsection}%
 {\normalfont\normalsize\scshape\noindent\thesubsection%
 \ifthenelse{\equal{#1}{}}%
 {}%
 {\ \letterspace{#1.}}%
 \ ---}%
}
\newcommand{\readme}{\par\vskip\baselineskip
 {\normalfont\normalsize\scshape\noindent%
  \letterspace{Readme.}\ ---}
}
\renewcommand\tableofcontents{%
 \sectionstar{\contentsname}%
 \@starttoc{toc}%
}
\renewcommand*\l@part[2]{%
 \addvspace{15pt \@plus\p@}%
 \noindent{\leavevmode%
  \scshape\letterspace{#1\qquad#2}%
 }\par\nobreak%
}
\renewcommand*\l@section[2]{%
 \setlength\@tempdima{\parindent}%
 \noindent
 {\leavevmode%
  \hskip\parindent#1\qquad#2%
 }\par\nobreak%
}
\makeatother
\usepackage{amsmath,amssymb}  
\usepackage{mathrsfs}         
\usepackage{mathabx}

\usepackage[thmmarks,amsmath]{ntheorem}
\usepackage{thmtools}

\numberwithin{equation}{subsection}

\def\sclsword#1{{\scshape\letterspace{#1}}}

\declaretheoremstyle[headformat=swapnumber,headpunct={.\ ---},%
headfont=\normalfont,bodyfont=\itshape,%
spaceabove=0pt,spacebelow=0pt,%
preheadhook={\vspace{15pt}},postfoothook={\vspace{-2pt}}]{theorem}
\declaretheorem[style=theorem,sibling=subsection,%
name=\sclsword{Theorem}]{theorem}
\declaretheorem[style=theorem,sibling=subsection,
name=\sclsword{Proposition}]{proposition}
\declaretheorem[style=theorem,sibling=subsection,
name=\sclsword{Lemma}]{lemma}
\declaretheorem[style=theorem,sibling=subsection,
name=\sclsword{Corollary}]{corollary}
\declaretheorem[style=theorem,sibling=subsection,
name=\sclsword{Conjecture}]{conjecture}

\declaretheoremstyle[headformat=swapnumber,headpunct={.\ ---},%
headfont=\normalfont,bodyfont=\normalfont,%
spaceabove=0pt,spacebelow=0pt,%
preheadhook={\vspace{15pt}},postfoothook={\vspace{-2pt}}]{definition}
\declaretheorem[style=definition,sibling=subsection,
name=\sclsword{Definition}]{definition}

\declaretheorem[style=definition,sibling=subsection,
name=\sclsword{Remark}]{remark}
\declaretheorem[style=definition,sibling=subsection,
name=\sclsword{Construction}]{construction}

\declaretheoremstyle[headpunct={\kern-5pt.},headfont=\normalfont,%
bodyfont=\normalfont,%
qed=\ensuremath{\square},spaceabove=0pt,spacebelow=0pt]{proof}
\declaretheoremstyle[headpunct={\kern-5pt.},headfont=\normalfont,%
bodyfont=\normalfont,%
qed=\ensuremath{\square},spaceabove=0pt,spacebelow=0pt]{nonumberproof}
\declaretheorem[style=proof,numbered=no,name=\emph{Proof}]{proof}

\declaretheoremstyle[headformat=swapnumber,headpunct={.\ ---},%
headfont=\normalfont,bodyfont=\normalfont,qed=\ensuremath{\square},%
spaceabove=0pt,spacebelow=0pt,%
preheadhook={\vspace{15pt}},postfoothook={\vspace{-2pt}}]{nproof}
\declaretheorem[style=nproof,sibling=subsection,name=\emph{Proof}]{nproof}

\usepackage{cleveref}
\crefname{condition}{condition}{conditions}
\crefname{conjecture}{conjecture}{conjectures}
\crefname{construction}{construction}{constructions}
\crefname{corollary}{corollary}{corollaries}
\crefname{diagram}{diagram}{diagrams}
\crefformat{subsection}{\S#2#1#3}
\crefformat{enumi}{\S#2#1#3}
\crefformat{nproof}{\S#2#1#3}
\creflabelformat{equation}{#2#1#3}

\newcommand{\into}{\hookrightarrow}
\newcommand{\onto}{\twoheadrightarrow}
\newcommand{\longto}{\longrightarrow}

\renewcommand{\Im}{\textnormal{Im}}

\newcommand{\Hom}{\textnormal{Hom}}
\newcommand{\End}{\textnormal{End}}

\newcommand{\iHom}{\underline{\Hom}}

\newcommand{\Mat}{\textnormal{Mat}}

\newcommand{\dual}[1]{\check{#1}}

\newcommand{\ZZ}{\mathbb{Z}}
\newcommand{\QQ}{\mathbb{Q}}
\newcommand{\QQl}{\QQ_{\ell}}

\newcommand{\QQp}{\QQ_{p}}

\newcommand{\RR}{\mathbb{R}}
\newcommand{\CC}{\mathbb{C}}
\newcommand{\HQ}{\mathbb{H}}
\newcommand{\FF}{\mathbb{F}}
\newcommand{\FFp}{\FF_{p}}
\newcommand{\FFq}{\FF_{q}}
\newcommand{\FFqbar}{\bar{\FF}_{q}}
\newcommand{\Adele}{\mathbb{A}}
\newcommand{\fin}{\textnormal{f}}

\newcommand{\primes}{\mathscr{L}}

\newcommand{\Spec}{\textnormal{Spec}}

\newcommand{\DelS}{\mathbb{S}}
\newcommand{\Sh}{\textnormal{Sh}}
\newcommand{\mSh}{\mathscr{S}}

\newcommand{\Kmpt}{\mathcal{K}}
\newcommand{\Sds}{\mathfrak{H}^{\pm}}
\newcommand{\AV}{\mathscr{A}}

\newcommand{\Gal}{\textnormal{Gal}}

\newcommand{\Vect}{\textnormal{Vect}}

\newcommand{\QHS}{\QQ\textnormal{HS}}

\makeatletter
\def\cpwith[#1]#2{\textnormal{c.p.}_{#1}(#2)}
\def\cpwithout#1{\textnormal{c.p.}(#1)}
\def\cp{\@ifnextchar[{\cpwith}{\cpwithout}}
\makeatother
\makeatletter
\def\Gmwith[#1]{\mathbb{G}_{\textnormal{m},#1}}
\def\Gmwithout{\mathbb{G}_{\textnormal{m}}}
\def\Gm{\@ifnextchar[{\Gmwith}{\Gmwithout}}
\makeatother
\newcommand{\GL}{\textnormal{GL}}
\newcommand{\GSp}{\textnormal{GSp}}

\newcommand{\ad}{\textnormal{ad}}

\newcommand{\cl}{\textnormal{cl}}

\newcommand{\et}{\textnormal{\'{e}t}}
\newcommand{\sing}{\textnormal{sing}}

\newcommand{\HH}{\textnormal{H}}
\newcommand{\Hl}{\HH_{\ell}}

\newcommand{\Hlambda}{\HH_{\lambda}}
\newcommand{\HB}{\HH_{\textnormal{B}}}

\newcommand{\Mot}{\textnormal{Mot}}

\newcommand{\GG}{\textnormal{G}}
\newcommand{\GB}{\GG_{\textnormal{B}}}

\newcommand{\Gl}{\GG_{\ell}}
\newcommand{\Glc}{\Gl^{\circ}}
\newcommand{\Glambda}{\GG_{\lambda}}
\newcommand{\ZB}{\textnormal{Z}_{\textnormal{B}}}
\newcommand{\Zl}{\textnormal{Z}_{\ell}}

\newcommand{\MTC}{\textnormal{MTC}}

\newcommand{\Res}{\textnormal{Res}}
\newcommand{\Nm}{\textnormal{Nm}}
\newcommand{\trace}{\textnormal{tr}}

\newcommand{\res}{\textnormal{res}}

\newcommand{\chrc}{\textnormal{char}}

\newcommand{\Tangen}[1]{\langle #1 \rangle^{\otimes}}

\newcommand{\tr}{\textsc{tr}}
\newcommand{\cm}{\textsc{cm}}
\usepackage{tikz}
\usetikzlibrary{cd,decorations.pathmorphing}

\def\title{On compatibility of the \(\ell\)-adic realisations\\
 of an abelian motive}
\def\author{Johan Commelin}

\usepackage{datetime}
\def\date{\dayofweekname{\day}{\month}{\year},
 the \ordinaldate{\day} of \monthname, \number\year}
\begin{document}
\begin{center}\Large\scshape
\letterspace{\title}
\end{center}

\medskip

\noindent\textit{by} \quad \author \hfill \date

\vskip3\baselineskip


\section{Introduction} 

\paragraph{Main result} 
The main result of this article is \cref{abmotCSR}:

{\narrower\it\noindent
Let \(M\) be an abelian motive
over a finitely generated subfield~\(K \subset \CC\).
Let \(E\) be a subfield of~\(\End(M)\),
and let \(\Lambda\) be the set of finite places of~\(E\).
Then the system \(\HH_{\Lambda}(M)\) is
a quasi-compatible system of Galois representations.\par
}

\noindent
To understand this result we need to explain what we mean by:
\begin{enumerate*}[label=(\textit{\roman*})]
 \item the words `abelian motive';
 \item the notation \(\HH_{\Lambda}(M)\); and
 \item the words `quasi-compatible system of Galois representations'.
\end{enumerate*}

\paragraph{Abelian motives} 
In this text we use motives in the sense of Andr\'{e}~\cite{An95}.
Alternatively we could have used the notion of absolute Hodge cycles.
Let \(K \subset \CC\) be a finitely generated subfield of the complex numbers.
An \emph{abelian motive} over~\(K\) is a summand of (a Tate~twist of)
the motive of an abelian variety over~\(K\).
In practice this means that an abelian motive~\(M\)
is a package consisting of a Hodge structure~\(\HB(M)\)
and for each prime~\(\ell\) an \(\ell\)-adic Galois representation~\(\Hl(M)\),
that arise in a compatible way as summands of Tate twists
of the cohomology of an abelian variety.

\paragraph{\(\lambda\)-adic realisations 
                                   and the notation \(\HH_{\Lambda}(M)\)\relax}
Let \(M\) be an abelian motive
over a finitely generated subfield~\(K \subset \CC\).
Let \(E\) be a subfield of~\(\End(M)\),
and let \(\Lambda\) be the set of finite places of the number field~\(E\).
For each prime number~\(\ell\),
the field~\(E\) acts on the Galois representation~\(\Hl(M)\)
via \(E_{\ell} = E \otimes \QQl = \prod_{\lambda|\ell} E_{\lambda}\)
and accordingly we get a decomposition
\(\Hl(M) = \bigoplus_{\lambda|\ell} \Hlambda(M)\)
of Galois representations,
where \(\Hlambda(M) = \Hl(M) \otimes_{E_{\ell}} E_{\lambda}\).
We denote with \(\HH_{\Lambda}(M)\)
the system of \(\lambda\)-adic Galois representations~\(\Hlambda(M)\)
as \(\lambda\) runs through~\(\Lambda\).

\paragraph{Quasi-compatible systems of Galois representations} \label{csrintro} 
In \cref{csogr} we develop a variation on the concept of
compatible systems of Galois representations
that has its origins in the work of Serre~\cite{Alrec}.
Besides the original work of Serre, we draw inspiration from Ribet~\cite{Ribet}
Larsen--Pink~\cite{LP2}, and Chi~\cite{Ch92}.
The main feature of our variant is a certain robustness
with respect to extension of the base field. 

By this we mean the following.
For the purpose of this introduction, let \(K\) be a number field.
(In \cref{csogr} the field~\(K\) is allowed to be any
finitely generated field of characteristic~\(0\).)
Let \(\rho_{\ell} \colon \Gal(\bar{K}/K) \to \GL_{n}(\QQ_{\ell})\)
and \(\rho_{\ell'} \colon \Gal(\bar{K}/K) \to \GL_{n}(\QQ_{\ell'})\)
be two Galois representations.
Let \(v\) be a place of~\(K\) such that the residue characteristic of~\(v\)
is neither~\(\ell\) nor~\(\ell'\).
Let \(\bar{v}\) be an extension of \(v\) to \(\bar{K}\).
Assume that \(\rho_{\ell}\) and~\(\rho_{\ell'}\)
are unramified at~\(\bar{v}/v\).
Let \(F_{\bar{v}/v}\) be a Frobenius element with respect to~\(v\)
(that is,
a lift of the Frobenius endomorphism of the residue field of~\(\bar{v}\)
to the decomposition group \(D_{\bar{v}/v} \subset \Gal(\bar{K}/K)\)).

We can now contrast the usual compatibility condition with our condition.
Recall that \(\rho_{\ell}\) and~\(\rho_{\ell'}\)
are called compatible at~\(v\)
if the characteristic polynomials of \(\rho_{\ell}(F_{\bar{v}/v})\)
and \(\rho_{\ell'}(F_{\bar{v}/v})\) have coefficients in~\(\QQ\)
and are equal to each other.
(Note that extensions of~\(v\) to~\(\bar{K}\) are conjugate to each other.
Consequently,
neither the condition that the representations are unramified
nor the compatibility condition on the characteristic polynomials
depends on the choice of~\(\bar{v}\).)

Our variant replaces the compatibility condition on the characteristic
polynomials of~\(\rho_{\ell}(F_{\bar{v}/v})\)
and~\(\rho_{\ell'}(F_{\bar{v}/v})\)
by the analogous condition for a power of~\(F_{\bar{v}/v}\)
that is allowed to depend on~\(v\).
In other words, one may replace the decomposition group~\(D_{\bar{v}/v}\)
by a finite index subgroup;
and in yet other words, one may replace the local field~\(K_{v}\)
by a finite field extension before checking the compatibility.
If \(\rho_{\ell}\) and~\(\rho_{\ell'}\) satisfy this relaxed condition,
then we say that they are \emph{quasi-compatible at~\(v\)}.

We may also take endomorphisms into account.
Instead of only considering systems of Galois representations
that are indexed by finite places of~\(\QQ\),
we may consider systems
that are indexed by finite places of a number field~\(E\).
This was already suggested by Serre~\cite{Alrec},
and Ribet pursued this further in~\cite{Ribet}.

The quasi-compatibility condition mentioned in the previous item
must then be adapted as follows.
Let \(\rho_{\lambda} \colon \Gal(\bar{K}/K) \to \GL_{n}(E_{\lambda})\)
and \(\rho_{\lambda'} \colon \Gal(\bar{K}/K) \to \GL_{n}(E_{\lambda'})\)
be two Galois representations.
Assume that the residue characteristic of~\(v\)
is different from the residue characteristics of~\(\lambda\) and~\(\lambda'\),
and assume that \(\rho_{\lambda}\) and~\(\rho_{\lambda'}\)
are unramified at~\(v\).
We say that \(\rho_{\lambda}\) and~\(\rho_{\lambda'}\)
are quasi-compatible at~\(v\)
if there is a positive integer~\(n\)
such that the characteristic polynomials of
\(\rho_{\lambda}(F_{\bar{v}/v}^{n})\)
and
\(\rho_{\lambda'}(F_{\bar{v}/v}^{n})\)
have coefficients in~\(E\) and are equal to each other.

The system \(\HH_{\Lambda}(M)\) mentioned above
is \emph{quasi-compatible} if for all \(\lambda,\lambda' \in \Lambda\)
there is a non-empty Zariski open subset \(U \subset \Spec(\mathcal{O}_{K})\)
such that \(\HH_{\lambda}(M)\) and \(\HH_{\lambda'}(M)\)
are quasi-compatible at all places \(v \in U\).

\paragraph{Related work} 
In~\cite{Lask14}, Laskar obtained related results.
Though he does not state this explicitly,
his results imply that for a large class of abelian motives
the \(\lambda\)-adic realisations form a
compatible system of Galois representations in the sense of Serre.
The contribution of the main result in this paper
is that \(M\) may be an arbitrary abelian motive;
although we need to weaken the concept of compatibility
to quasi-compatibility to achieve this.
See~\cref{compLask} for more details.

\paragraph{Organisation of the paper} 
Every section starts with a paragraph labeled `\textsc{\letterspace{Readme}}'.
These paragraphs highlight the important parts of their section,
or describe the role of the section in the text as a whole.
We hope that these paragraphs aid in navigating the text.

In \cref{abmot} we recall the definition of abelian motives
in the sense of Andr\'{e}~\cite{An95}.
We also recall useful properties of abelian motives.
This section does not contain new results.
In \cref{csogr} we give the main definition of this paper,
namely the notion of a quasi-compatible system of Galois representations.
In \cref{CSReg} we recall results showing that abelian varieties and
so-called \cm~motives give rise to such quasi-compatible systems.
These results are known over number fields,
and we make the rather trivial generalisation to finitely generated fields.
\Cref{mainproof} is the heart of this paper,
as it proves the main result.
See below for an outline of its contents.
Finally, \cref{props} and \cref{MTCrmk} are appendices.
In the former we show that quasi-compatible systems
share some of the familiar properties of
compatible systems in the sense of Serre.
The latter appendix shows that
for abelian motives the Mumford--Tate conjecture
does not depend on the prime number~\(\ell\) that occurs in its statement.

\paragraph{Outline of the proof} 
Shimura showed that if \(M = \HH^{1}(A)\),
with \(A\) an abelian variety,
then the system \(\HH_{\Lambda}(M)\) is
an \(E\)-rational compatible system in the sense of Serre.
In \cref{abvarCSR} we recall this result of Shimura
in the setting of quasi-compatible systems of Galois representations.

\Cref{mainproof} proves the main result of this paper,
namely that \(\HH_{\Lambda}(M)\) is
a quasi-compatible system of Galois representations
for every abelian motive~\(M\).
Roughly speaking, the proof works by placing the abelian motive~\(M\)
in a family of motives over a Shimura variety.
The problem may then be deformed to a \cm~point on the Shimura variety,
where we can prove the result by reducing to the case
of abelian varieties mentioned above.
To make this work we need a recent result of Kisin~\cite{Kisin_modp}:
Let \(\mSh\)~be an integral model of a Shimura variety of Hodge type
over the ring of integers~\(\mathcal{O}_{K}\) of a \(p\)-adic field~\(K\),
satisfying some additional technical conditions.
Then every point in the special fibre of~\(\mSh\) is \emph{isogenous} to
a point that lifts to a \cm~point of the generic fibre~\(\mSh_{K}\).
We refer to the main text for details (\cref{Kisin_isog} and~\cref{isoglift}).

\paragraph{Terminology and notation} \label{notation} 
We say that a field is a \emph{finitely generated field}
if it is finitely generated over its prime field.
A motive \(M\) over a field \(K \subset \CC\) is called
\emph{geometrically irreducible} if \(M_{\CC}\) is irreducible.
If \(G\) is a semiabelian variety,
then we denote with \(\End^{0}(G)\) the \(\QQ\)-algebra \(\End(G) \otimes \QQ\).
If \(X\) is a scheme,
then \(X^{\cl}\) denotes the set of closed points of~\(X\).

If \(K\) is a field, \(V\) a vector space over~\(K\),
and \(g\) an endomorphism of~\(V\),
then we denote with \(\cp[K]{g|V}\) the characteristic polynomial of~\(g\).
If there is no confusion possible,
then we may drop~\(K\) or~\(V\) from the notation,
and write \(\cp{g|V}\) or simply~\(\cp{g}\).

Let \(E\) be a number field.
Recall that \(E\) is called \emph{totally real}~(\tr)
if for all complex embeddings \(\sigma \colon E \into \CC\)
the image \(\sigma(E)\) is contained in~\(\RR\).
The field~\(E\) is called a \emph{complex multiplication} field~(\cm)
if it is a quadratic extension of a totally real field
(typically denoted~\(E^{0}\)),
and if all complex embeddings \(\sigma \colon E \into \CC\)
have an image that is not contained in~\(\RR\).

Let \(C\) be a Tannakian category, and let \(V\) be an object of~\(C\).
If \(a\) and~\(b\) are non-negative integers, then \(\mathbb{T}^{a,b}V\)
denotes the object \(V^{\otimes a} \otimes \dual{V}^{\otimes b}\).
With \(\Tangen{V}\) we denote the
smallest full Tannakian subcategory of~\(C\) that contains~\(V\).
This means that it is the smallest full subcategory of~\(C\)
that contains~\(V\) and that is closed under
directs sums, tensor products, duals, and subquotients.
The irreducible objects in \(\Tangen{V}\) are precisely
the irreducible objects of~\(C\) that are a subquotient
of~\(\mathbb{T}^{a,b}V\) for some \(a,b \ge 0\).

\paragraph{Acknowledgements} 
This paper is part of the author's PhD~thesis~\cite{Comm2017PhD}.
I warmly thank Ben Moonen for his thorough and stimulating supervision.
I also thank him for pointing me to Kisin's paper~\cite{Kisin_modp}.
I thank Netan Dogra and Milan Lopuha\"{a} for several useful discussions.
Anna Cadoret, Bas Edixhoven, Rutger Noot, and Lenny Taelman
have provided extensive feedback, for which I thank them.
In particular, Rutger Noot pointed me to the work of Laskar~\cite{Lask14}.
I express gratitude to Pierre Deligne for his interest in this work;
his valuable comments have improved several parts of this paper.

This research has been financially supported by the Netherlands Organisation
for Scientific Research~(NWO) under project no.~613.001.207
\emph{(Arithmetic and motivic aspects of the Kuga--Satake construction)}.

\section{Abelian motives} \label{abmot} 

\readme 
We briefly review the definition of abelian motives
in the sense of Andr\'{e}~\cite{An95},
and we recall some of their useful properties.

\paragraph{} \label{Weil_cohoms} 
Let \(K \subset \CC\) be a field.
Let \(X\) be a smooth projective variety over~\(K\).
For every prime number~\(\ell\),
let \(\Hl^{i}(X)\) denote the Galois representation
\(\HH_{\et}^{i}(X_{\bar{K}}, \QQl)\).
Write \(\HB^{i}(X)\) for the Hodge structure \(\HH_{\sing}^{i}(X(\CC), \QQ)\).

\paragraph{} 
Let \(K \subset \CC\) be a field.
In this text a \emph{motive} over~\(K\) shall mean
a motive in the sense of Andr\'{e}~\cite{An95}.
(To be precise, our category of base pieces is the category of
smooth projective varieties over~\(K\),
and our reference cohomology is Betti cohomology, \(\HB(\_)\).
The resulting notion of motive does not depend on the chosen
reference cohomology, see proposition~2.3 of~\cite{An95}.)
We denote the category of motives over~\(K\) with \(\Mot_{K}\).

The category~\(\Mot_{K}\) is a \emph{semisimple} neutral Tannakian category
and therefore the motivic Galois group of a motive
is a reductive algebraic group.
We further mention that K\"{u}nneth projectors exist in~\(\Mot_{K}\).
If \(K = \CC\), then we know that the Betti realisation functor is fully
faithful on the Tannakian subcategory generated by motives of abelian
varieties, see~\cref{hdgmot}.

\paragraph{} \label{realisations} 
Let \(K \subset \CC\) be a field.
If \(X\) is a smooth projective variety over~\(K\),
then we write \(\HH^{i}(X)\) for the motive in degree~\(i\)
associated with~\(X\).
The cohomology functors mentioned in~\cref{Weil_cohoms}
induce realisation functors on the category of motives over~\(K\).
Let \(M\) be a motive over~\(K\).
For every prime~\(\ell\),
we write \(\Hl(M)\) for the \(\ell\)-adic realisation;
it is a finite-dimensional \(\QQl\)-vector space
equipped with a continuous representation of~\(\Gal(\bar{K}/K)\).
Similarly, we write \(\HB(M)\) for the Betti realisation;
it is a polarisable \(\QQ\)-Hodge structure.

\paragraph{} \label{Artin} 
Let \(M\) be a motive over~\(\CC\), and let \(\ell\) be a prime number.
There is an isomorphism of \(\QQl\)-vector spaces
\(\HB(M) \otimes_{\QQ} \QQl \cong \Hl(M)\).
If \(K \subset \CC\) is a subfield such that \(M\) is defined over~\(K\),
then there is an isomorphism of \(\QQl\)-vector spaces
\(\Hl(M) \cong \Hl(M_{\CC})\);
and therefore
\[
 \HB(M) \otimes_{\QQ} \QQl \cong \Hl(M).
\]
This isomorphism was proven for varieties
by Artin in expos\'{e}~\textsc{xi} in~\cite{SGA4-3}.
The generalisation to motives follows from the fact
that the isomorphism is compatible with cycle class maps.

\paragraph{} \label{MTgrp} 
Let \(V\) be a \(\QQ\)-Hodge structure.
The \emph{Mumford--Tate group}~\(\GB(V)\) of~\(V\)
is the linear algebraic group over~\(\QQ\) associated with the
Tannakian category \(\Tangen{V}\) generated by~\(V\)
(with the forgetful functor \(\QHS \to \Vect_{\QQ}\) as fibre functor).
If \(V\) is polarisable,
then the Tannakian category \(\Tangen{V}\) is semisimple;
which implies that \(\GB(V)\) is reductive.

For an alternative description, recall that the Hodge structure on~\(V\)
is determined by a homomorphism of algebraic groups
\(\DelS \to \GL(V \otimes_{\QQ} \RR)\),
where \(\DelS\) is the Deligne torus \(\Res^{\CC}_{\RR}\Gm\).
The Mumford--Tate group is the smallest algebraic subgroup \(G\) of \(\GL(V)\)
such that \(G_{\RR}\) contains the image of~\(\DelS\).
Since \(\DelS\) is connected, so is~\(\GB(V)\).

If \(M\)~is a motive over a field \(K \subset \CC\),
then we write \(\GB(M)\) for \(\GB(\HB(M))\).

\paragraph{} \label{abmotdef} 
Let \(K \subset \CC\) be a field.
An \emph{abelian motive} over~\(K\) is an object of the Tannakian subcategory of
motives over~\(K\) generated by the motives of abelian varieties over~\(K\).
Recall that \(\HH(A) \cong \bigwedge^{\star} \HH^{1}(A)\)
for every abelian variety~\(A\) over~\(K\),
and thus we have \(\Tangen{\HH(A)} = \Tangen{\HH^{1}(A)}\).
If \(A\) is a non-trivial abelian variety,
then the class of any effective non-zero divisor
realises \(\mathds{1}(-1)\) as a subobject of~\(\HH^{2}(A)\),
and therefore \(\mathds{1}(-1) \in \Tangen{\HH^{1}(A)}\).
In particular \(\mathds{1}(-1)\) is an abelian motive.
We claim that every abelian motive~\(M\) is
contained in \(\Tangen{\HH^{1}(A)}\) for some abelian variety~\(A\) over~\(K\).
By definition there are abelian varieties \((A_{i})_{i = 1}^{k}\)
such that \(M\) is contained in the Tannakian subcategory
generated by the \(\HH(A_{i})\).
Put \(A = \prod_{i = 1}^{k} A_{i}\),
so that \(\HH^{1}(A) \cong \bigoplus_{i = 1}^{k} \HH^{1}(A_{i})\).
It follows that \(M\) is contained in \(\Tangen{\HH^{1}(A)}\).

\begin{theorem} \label{hdgmot} 
  The Betti realisation functor \(\HB(\_)\) is fully faithful on the
  subcategory of abelian motives over~\(\CC\).
 \begin{proof}
   See th\'{e}or\`{e}me~0.6.2 of~\cite{An95}.
 \end{proof}
\end{theorem}

\paragraph{} 
In the rest of this section we focus on so-called \cm~motives.
They will play a crucial r\^{o}le in the proof of our main result.
An important tool in understanding abelian \cm~motives
is the half-twist construction that we describe in \cref{halftwists}.

\begin{definition} 
 A motive~\(M\) over a field~\(K \subset \CC\)
 is called a \emph{\cm~motive}
 if \(\HB(M)\) is a \cm~Hodge structure
 (\textit{i.e.}, the group \(\GB(M)\) is commutative).
\end{definition}

\paragraph{} 
Let \(E\) be a \cm~field.
Let \(\Sigma(E)\) be the set of complex embeddings of~\(E\).
The complex conjugation on \(E\) induces
an involution \(\sigma \mapsto \sigma^{\dagger}\) on~\(\Sigma(E)\).
If \(T\) is a subset of~\(\Sigma(E)\),
then we denote with~\(T^{\dagger}\) the image of~\(T\) under this involution.
Recall that a \cm~type \(\Phi \subset \Sigma(E)\) is a subset
such that \(\Phi \cup \Phi^{\dagger} = \Sigma(E)\)
and \(\Phi \cap \Phi^{\dagger} = \varnothing\).
Each \cm~type \(\Phi\) defines
a Hodge structure \(E_{\Phi}\) on~\(E\) of type \(\{(0,1), (1,0)\}\),
via
\[
 E_{\Phi} \otimes_{\QQ} \CC \cong \CC^{\Sigma(E)},
 \qquad E_{\Phi}^{0,1} \cong \CC^{\Phi^{\dagger}},
 \qquad E_{\Phi}^{1,0} \cong \CC^{\Phi}.
\]

\paragraph{Half-twists} \label{halftwists} 
The idea of half-twists originates from~\cite{vG01},
though we use the description in~\S7 of~\cite{Mo15}.
Let \(V\) be a Hodge structure of weight~\(n\).
The level of~\(V\), denoted~\(m\),
is by definition \(\max\{ p-q \mid V^{p,q} \ne 0\}\).
Suppose that \(\End(V)\) contains a \cm~field~\(E\).
Let \(\Sigma(E)\) denote the set of complex embeddings \(E \into \CC\).
Let \(T \subset \Sigma(E)\) be the embeddings
through which \(E\) acts on \(\bigoplus_{p \ge \lceil n/2 \rceil} V^{p,q}\).
Assume that \(T \cap T^{\dagger} = \varnothing\).
(Note that if \(\dim_{E}(V) = 1\),
then the condition \(T \cap T^{\dagger} = \varnothing\) is certainly satisfied.)

Let \(\Phi \subset \Sigma(E)\) be a \cm~type,
and let \(E_{\Phi}\) be the associated Hodge structure on~\(E\).
If \(T \cap \Phi = \varnothing\) and \(m \ge 1\),
then the Hodge structure \(W = E_{\Phi} \otimes_{E} V\)
has weight~\(n+1\) and level~\(m-1\).
In that case we call \(W\) a \emph{half-twist} of~\(V\).
Note that under our assumption \(T \cap T^{\dagger} = \varnothing\)
we can certainly find a \cm~type with \(T \cap \Phi = \varnothing\),
so that there exist half-twists of~\(V\).
For each \cm~type~\(\Phi\) with \(T \cap \Phi = \varnothing\),
there is a complex abelian variety~\(A_{\Phi}\)
(well-defined up to isogeny),
with \(\HB^{1}(A_{\Phi}) \cong E_{\Phi}\).
By construction we have \(E \subset \End(\HB^{1}(A_{\Phi}))\)
and \(E \subset \End(W)\).
Note that \(V \cong \iHom_{E}(\HB^{1}(A_{\Phi}), W)\).
In the next paragraph we will see
that this construction generalises to abelian motives.

\paragraph{} \label{halftwistsMot} 
Let \(K \subset \CC\) be a field.
Let \(M\) be an abelian motive over~\(K\).
Assume that \(M\) is pure of weight~\(n\),
and assume that \(\End(M)\) contains a \cm~field~\(E\).
Note that \(\HB(M)\) is a Hodge structure of weight~\(n\).
Let \(T \subset \Sigma(E)\) be the set of embeddings through which \(E\) acts
on~\(\bigoplus_{p \ge \lceil n/2 \rceil}\HB(M)^{p,q}\).
Assume that \(T \cap T^{\dagger} = \varnothing\).
Then there exists a finitely generated extension \(L/K\),
an abelian variety \(A\) over~\(L\), and a motive~\(N\) over~\(L\),
such that \(E \subset \End(\HH^{1}(A))\), and \(E \subset \End(N)\),
and such that \(M_{L} \cong \iHom_{E}(\HH^{1}(A), N)\).
Indeed, choose a \cm~type \(\Phi \subset \Sigma(E)\)
such that \(T \cap \Phi = \varnothing\).
Put \(A = A_{\Phi}\), and \(N = \HH^{1}(A_{\Phi}) \otimes_{E} M_{\CC}\).
Then \(M_{\CC} \cong \iHom_{E}(\HH^{1}(A), N)\),
by \cref{hdgmot} and the construction above.
The abelian variety~\(A\) and the motive~\(N\) are defined
over some finitely generated extension of~\(K\), and so is the isomorphism
\(M_{\CC} \cong \iHom_{E}(\HH^{1}(A), N)\).

\section{Quasi-compatible systems of Galois representations} 
\label{csogr}

\readme 
In this section we develop the notion of
quasi-compatible systems of Galois representations,
a variant on Serre's compatible systems of Galois representations~\cite{Alrec}.
We follow Serre's suggestion of developing
an \(E\)-rational version (where \(E\) is a number field);
which has also been done by Ribet~\cite{Ribet} and Chi~\cite{Ch92}.
The main benefit of the variant that we develop is
that we relax the compatibility condition,
thereby gaining a certain robustness
with respect to extensions of the base field and residue fields.
We will need this property in a crucial way in the proof of \cref{abmotCSR}.

\paragraph{} 
Let \(\kappa\) be a finite field with~\(q\) elements,
and let \(\bar{\kappa}\) be an algebraic closure of~\(\kappa\).
We denote with \(F_{\bar{\kappa}/\kappa}\) the
geometric Frobenius element
(that is, the inverse of \(x \mapsto x^{q}\))
in \(\Gal(\bar{\kappa}/\kappa)\).

\paragraph{} \label{numfldFrob} 
Let \(K\) be a number field.
Let \(v\) be a finite place of~\(K\),
and let \(K_{v}\) denote the completion of~\(K\) at~\(v\).
Let \(\bar{K}_{v}\) be an algebraic closure of~\(K_{v}\).
Let \(\bar{\kappa}/\kappa\) be the extension of residue fields
corresponding with \(\bar{K}_{v}/K_{v}\).
The inertia group, denoted~\(I_{v}\),
is the kernel of the natural surjection
\(\Gal(\bar{K}_{v}/K_{v}) \onto \Gal(\bar{\kappa}/\kappa)\).
The inverse image of \(F_{\bar{\kappa}/\kappa}\)
in \(\Gal(\bar{K}_{v}/K_{v})\)
is called the \emph{Frobenius coset} of~\(v\).
An element \(\alpha \in \Gal(\bar{K}/K)\) is called a
\emph{Frobenius element with respect to~\(v\)}
if there exists an embedding \(\bar{K} \into \bar{K}_{v}\)
such that \(\alpha\) is the restriction
of an element of the Frobenius coset of~\(v\).

\paragraph{} 
Let \(K\) be a finitely generated field.
A \emph{model} of~\(K\) is an
integral scheme~\(X\) of finite type over~\(\Spec(\ZZ)\)
together with an isomorphism between \(K\) and the function field of~\(X\).
Remark that if \(K\) is a number field,
and \(R \subset K\) is an order,
then \(\Spec(R)\) is naturally a model of~\(K\).
The only model of a number field~\(K\) that is normal
and proper over~\(\Spec(\ZZ)\) is \(\Spec(\mathcal{O}_{K})\).

\paragraph{} \label{Frobxmodel} 
Let \(K\) be a finitely generated field,
and let \(X\) be a model of~\(K\).
Recall that we denote the set of closed points of~\(X\) with~\(X^{\cl}\).
Let \(x \in X^{\cl}\) be a closed point.
Let \(K_{x}\) be the function field of the Henselisation of \(X\) at~\(x\);
and let \(\kappa(x)\) be the residue field at~\(x\).
We denote with \(I_{x}\) the kernel of
\(\Gal(\bar{K}_{x}/K_{x}) \onto \Gal(\bar{\kappa}(x)/\kappa(x))\).
Every embedding \(\bar{K} \into \bar{K}_{x}\)
induces an inclusion \(\Gal(\bar{K}_{x}/K_{x}) \into \Gal(\bar{K}/K)\).

Like in \cref{numfldFrob},
the inverse image of \(F_{\bar{\kappa}(x)/\kappa(x)}\)
in \(\Gal(\bar{K}_{x}/K_{x})\)
is called the Frobenius coset of~\(x\).
An element \(\alpha \in \Gal(\bar{K}/K)\) is called a
\emph{Frobenius element with respect to~\(x\)}
if there exists an embedding \(\bar{K} \into \bar{K}_{x}\)
such that \(\alpha\) is the restriction
of an element of the Frobenius coset of~\(x\).

\begin{definition} 
 Let \(K\) be a field, let \(E\) be a number field,
 and let \(\lambda\) be a place of~\(E\).
 A \emph{\(\lambda\)-adic Galois representation} of~\(K\)
 is a representation of \(\Gal(\bar{K}/K)\)
 on a finite-dimensional \(E_{\lambda}\)-vector space
 that is continuous for the \(\lambda\)-adic topology.

 Let \(\rho_{\lambda} \colon \Gal(\bar{K}/K) \to \GL(V_{\lambda})\)
 be a \(\lambda\)-adic Galois representation of~\(K\).
 We denote with \(\Glambda(\rho_{\lambda})\) or \(\Glambda(V_{\lambda})\) the
 Zariski closure of the image of \(\Gal(\bar{K}/K)\) in \(\GL(V_{\lambda})\).
 In particular, if \(E = \QQ\) and \(\lambda = \ell\),
 then we denote this group with \(\Gl(\rho_{\ell})\) or~\(\Gl(V_{\ell})\).
\end{definition}

\paragraph{} 
Let \(K\) be a finitely generated field.
Let \(X\) be a model of~\(K\), and let \(x \in X^{\cl}\) be a closed point.
We use the notation introduced in \cref{Frobxmodel}.
Let \(\rho\) be a \(\lambda\)-adic Galois representation of~\(K\).
We say that \(\rho\) is \emph{unramified at~\(x\)}
if there is an embedding \(\bar{K} \into \bar{K}_{x}\)
for which \(\rho(I_{x}) = \{1\}\).
If this is true for one embedding, then it is true for all embeddings.

Let \(F_{x}\) be a Frobenius element with respect to~\(x\).
If \(\rho\) is unramified at~\(x\),
then the element \(F_{x,\rho} = \rho(F_{x})\) is well-defined up to conjugation.
For \(n \in \ZZ\), we write \(P_{x,\rho,n}(t)\)
for the characteristic polynomial \(\cp{F_{x,\rho}^{n}}\).
Note that \(P_{x,\rho,n}(t)\) is well-defined,
since conjugate endomorphisms have the same characteristic polynomial.

\paragraph{} 
In the following definitions,
one recovers the notions of Serre~\cite{Alrec}
by demanding \(n = 1\) everywhere.
By not making this demand we gain a certain flexibility
that will turn out to be crucial for our proof of \cref{abmotCSR}.

\begin{definition} \label{Erat} 
 Let \(K\) be a finitely generated field.
 Let \(E\) be a number field,
 and let \(\lambda\) be a finite place of~\(E\).
 Let \(\rho\) be a \(\lambda\)-adic Galois representation of~\(K\).
 Let \(X\) be a model of~\(K\),
 and let \(x \in X^{\cl}\) be a closed point.
 The representation~\(\rho\) is said to be \emph{\(E\)-rational at~\(x\)} if
 \(\rho\)~is unramified at~\(x\), and \(P_{x,\rho,n}(t) \in E[t]\),
 for some \(n \ge 1\).
\end{definition}

\begin{definition} \label{comprep} 
 Let \(K\) be a finitely generated field.
 Let \(E\) be a number field, and
 let \(\lambda_{1}\) and~\(\lambda_{2}\) be two finite places of~\(E\).
 For \(i = 1,2\), let \(\rho_{i}\) be
 a \(\lambda_{i}\)-adic Galois representation of~\(K\).
 \begin{enumerate}
  \item Let \(X\) be a model of~\(K\),
   and let \(x \in X^{\cl}\) be a closed point.
   Then \(\rho_{1}\) and~\(\rho_{2}\) are said to be
   \emph{quasi-compatible at~\(x\)}
   if \(\rho_{1}\) and~\(\rho_{2}\) are both \(E\)-rational at~\(x\),
   and if there is an integer \(n\) such that
   \(P_{x,\rho_{1},n}(t) = P_{x,\rho_{2},n}(t)\)
   as polynomials in \(E[t]\).
  \item \label[definition]{compatwrtmodel}
   Let \(X\) be a model of~\(K\).
   The representations \(\rho_{1}\) and~\(\rho_{2}\) are
   \emph{quasi-compatible with respect to~\(X\)}
   if there is a non-empty open subset \(U \subset X\),
   such that \(\rho_{1}\) and~\(\rho_{2}\) are quasi-compatible at~\(x\)
   for all \(x \in U^{\cl}\).
  \item \label[definition]{strongcompatwrtmodel}
   Let \(X\) be a model of~\(K\).
   The representations \(\rho_{1}\) and~\(\rho_{2}\) are
   \emph{strongly quasi-compatible with respect to~\(X\)}
   if \(\rho_{1}\) and~\(\rho_{2}\) are quasi-compatible
   at all points \(x \in X^{\cl}\) that satisfy the following condition:
   
   {\narrower\noindent
    The places \(\lambda_{1}\) and~\(\lambda_{2}\)
    have a residue characteristic that is different from the
    residue characteristic of~\(x\),
    and \(\rho_{1}\) and~\(\rho_{2}\)
    are unramified at~\(x\).
    \par}
  \item \label[definition]{compatallmodel}
   The representations \(\rho_{1}\) and~\(\rho_{2}\) are
   \emph{(strongly) quasi-compatible}
   if they are (strongly) quasi-compatible with respect to every model of~\(K\).
 \end{enumerate}
\end{definition}

\begin{remark} \label{compatbirat} 
 Let \(K\), \(E\), \(\lambda_{1}\), \(\lambda_{2}\),
 \(\rho_{1}\), and \(\rho_{2}\) be as in the above definition.
 \begin{enumerate}
  \item If there is one model~\(X\) of~\(K\)
   such that \(\rho_{1}\) and~\(\rho_{2}\)
   are quasi-compatible with respect to~\(X\),
   then \(\rho_{1}\) and~\(\rho_{2}\)
   are quasi-compatible with respect to every model of~\(K\),
   since all models of~\(K\) are birational to each other.
  \item It is \emph{not} known whether the notion of strong quasi-compatibility
   is stable under birational equivalence:
   if \(\rho_{1}\) and~\(\rho_{2}\) are quasi-compatible
   with respect to some model~\(X\) of~\(K\),
   then by definition there exists a non-empty open subset
   \(U \subset X\) such that \(\rho_{1}\) and~\(\rho_{2}\)
   are strongly quasi-compatible with respect to~\(U\).
   But there is no \emph{a priori} reason to expect that
   \(\rho_{1}\) and~\(\rho_{2}\)
   are strongly quasi-compatible with respect to~\(X\).
  \item It is \emph{not} known whether
   strong quasi-compatibility is an equivalence relation:
   Let \(\rho_{1}\),~\(\rho_{2}\), and~\(\rho_{3}\)
   be respectively \(\lambda_{1}\)-adic,
   \(\lambda_{2}\)-adic, and~\(\lambda_{3}\)-adic
   Galois representations of~\(K\).
   Suppose that \(\rho_{1}\) and~\(\rho_{2}\) are strongly compatible
   and
   suppose that \(\rho_{2}\) and~\(\rho_{3}\) are strongly compatible.
   Then it is \emph{not} known whether
   \(\rho_{1}\) and~\(\rho_{3}\) are strongly compatible.
 \end{enumerate}
\end{remark}

\begin{definition} 
 Let \(K\) be a field.
 With a \emph{system of Galois representations} of~\(K\) we mean a triple
 \((E, \Lambda, (\rho_{\lambda})_{\lambda \in \Lambda})\),
 where
 \(E\) is a number field;
 \(\Lambda\) is a set of finite places of~\(E\); and
 \(\rho_{\lambda}\) (\(\lambda \in \Lambda\))
 is a \(\lambda\)-adic Galois representation of~\(K\).
\end{definition}

\paragraph{} 
In what follows, we often denote a system of Galois representations
\((E, \Lambda, (\rho_{\lambda})_{\lambda \in \Lambda})\)
with \(\rho_{\Lambda}\), leaving the number field~\(E\) implicit.
In contexts where there are multiple number fields
the notation will make clear which number field is meant
(\textit{e.g.,} by denoting the set of finite places of a number field~\(E'\)
with~\(\Lambda'\), etc\ldots).

\begin{definition} \label{CSR} 
 Let \(K\) be a finitely generated field.
 Let \(E\) be a number field,
 and let \(\Lambda\) be a set of finite places of~\(E\).
 Let \(\rho_{\Lambda}\) be a system of Galois representations of~\(K\).
 \begin{enumerate}
  \item Let \(X\) be a model of~\(K\).
   The system~\(\rho_{\Lambda}\) is
   \emph{(strongly) quasi-compatible with respect to~\(X\)}
   if for all \(\lambda_{1}, \lambda_{2} \in \Lambda\)
   the representations \(\rho_{\lambda_{1}}\) and~\(\rho_{\lambda_{2}}\)
   are (strongly) quasi-compatible with respect to~\(X\).
  \item The system~\(\rho_{\Lambda}\) is called
   \emph{(strongly) quasi-compatible}
   if for all \(\lambda_{1}, \lambda_{2} \in \Lambda\)
   the representations \(\rho_{\lambda_{1}}\) and~\(\rho_{\lambda_{2}}\)
   are (strongly) quasi-compatible.
 \end{enumerate}
\end{definition}

\begin{remark} 
 The first two points of \cref{compatbirat}
 apply mutatis mutandis to the above definition:
 compatibility is stable under birational equivalence,
 but for strong compatibility we do not know this.
\end{remark}

\begin{lemma} \label{CSRbc} 
 Let \(K\) be a finitely generated field.
 Let \(E\) be a number field,
 and let \(\Lambda\) be a set of finite places of~\(E\).
 Let \(\rho_{\Lambda}\) be a system of Galois representations of~\(K\).
 Let \(L\) be a finitely generated extension of~\(K\).
 Let \(\rho'_{\Lambda}\) denote the system of Galois representations of~\(L\)
 obtained by restricting the system~\(\rho_{\Lambda}\) to~\(L\).
 \begin{enumerate}
  \item The system~\(\rho_{\Lambda}\) is quasi-compatible
   if and only if
   the system~\(\rho'_{\Lambda}\) is quasi-compatible.
  \item If the system~\(\rho'_{\Lambda}\) is strongly quasi-compatible,
   then the system~\(\rho_{\Lambda}\) is strongly quasi-compatible.
 \end{enumerate}
 \begin{proof}
  Without loss of generality we may and do assume that
  \(\Lambda = \{\lambda_{1}, \lambda_{2}\}\).
  Let \(X\) be a model of~\(K\).
  Let \(Y\) be an~\(X\)-scheme that is a model of~\(L\).
  Let \(x \in X^{\cl}\) be a closed point whose residue characteristic is
  different from the residue characteristic
  of~\(\lambda_{1}\) and~\(\lambda_{2}\).

  For the remainder of the proof,
  we may and do assume that \(\rho_{\lambda_{1}}\) and~\(\rho_{\lambda_{2}}\)
  are both unramified at~\(x\).
  Then \(\rho'_{\lambda_{1}}\) and~\(\rho'_{\lambda_{2}}\)
  are both unramified at all points~\(y \in Y_{x}^{\cl}\).
  If \(y \in Y_{x}^{\cl}\) is a closed point,
  and \(k\) denotes the residue extension degree \([\kappa(y) : \kappa(x)]\),
  then we have \(F_{y,\rho_{\lambda}} = F_{x,\rho_{\lambda}}^{k}\)
  for all \(\lambda \in \Lambda\).
  This leads to the following conclusions:
  \begin{enumerate*}[label=(\textit{\roman*})]
   \item For every point \(y \in Y_{x}^{\cl}\),
    if \(\rho'_{\lambda_{1}}\) and~\(\rho'_{\lambda_{2}}\)
    are quasi-compatible at~\(y\),
    then \(\rho_{\lambda_{1}}\) and~\(\rho_{\lambda_{2}}\)
    are quasi-compatible at~\(x\); and
   \item if \(\rho_{\lambda_{1}}\) and~\(\rho_{\lambda_{2}}\)
    are quasi-compatible at~\(x\),
    then \(\rho'_{\lambda_{1}}\) and~\(\rho'_{\lambda_{2}}\)
    are quasi-compatible at all points \(y \in Y_{x}^{\cl}\).
  \end{enumerate*}
  Together, these two conclusions complete the proof.
 \end{proof}
\end{lemma}
(Note that I cannot prove the converse implication in point~2,
for the following reason.
Let \(y \in Y^{\cl}\) be a closed point whose residue characteristic is
different from the residue characteristic
of~\(\lambda_{1}\) and~\(\lambda_{2}\).
If \(\rho'_{\lambda_{1}}\) and~\(\rho'_{\lambda_{2}}\)
are unramified at~\(y\),
but \(\rho_{\lambda_{1}}\) and~\(\rho_{\lambda_{2}}\)
are not unramified at the image~\(x\) of~\(y\) in~\(X\),
then I~do not see how to prove that
\(\rho'_{\lambda_{1}}\) and~\(\rho'_{\lambda_{2}}\)
are quasi-compatible at~\(y\).)

\paragraph{} \label{restrictE} 
Let \(K\) be a finitely generated field.
Let \(E\) be a number field, and
let \(\Lambda\) be a set of finite places of~\(E\).
Let \(\rho_{\Lambda}\) be a system of Galois representations over~\(K\).
Let \(E' \subset E\) be a subfield,
and let \(\Lambda'\) be the set of places \(\lambda'\) of~\(E'\)
satisfying the following condition:

{\narrower\noindent
 For all places \(\lambda\) of~\(E\) with \(\lambda | \lambda'\),
 we have \(\lambda \in \Lambda\).
 \par}

\noindent%
For each \(\lambda' \in \Lambda'\), the representation
\(\rho_{\lambda'} = \bigoplus_{\lambda | \lambda'} \rho_{\lambda}\)
is naturally a \(\lambda'\)-adic Galois represenation of~\(K\).
We thus obtain a system of Galois representations~\(\rho_{\Lambda'}\).

\begin{lemma} \label{restrictElem} 
 Let \(K\), \(E' \subset E\), \(\Lambda\), \(\Lambda'\),
 \(\rho_{\Lambda}\), and~\(\rho_{\Lambda'}\) be as in \cref{restrictE}.
 If \(\rho_{\Lambda}\) is
 a (strongly) quasi-compatible system of Galois representations,
 then \(\rho_{\Lambda'}\) is
 a (strongly) quasi-compatible system of Galois representations.
 \begin{proof}
  To see this, we may assume that \(\Lambda' = \{\lambda'_{1}, \lambda'_{2}\}\)
  and \(\Lambda\) is the set of all places~\(\lambda\) of~\(E\)
  that lie above a place \(\lambda' \in \Lambda'\).
  Let \(X\) be a model of~\(K\).
  Let \(x \in X^{\cl}\) be a closed point
  whose residue characteristic is different from
  the residue characteristic of~\(\lambda'_{1}\) and~\(\lambda'_{2}\).
  Assume that \(\rho_{\lambda'_{1}}\) and~\(\rho_{\lambda'_{2}}\)
  are both unramified at~\(x\).
  Suppose that for all \(\lambda_{1}, \lambda_{2} \in \Lambda\),
  the representations \(\rho_{\lambda_{1}}\) and~\(\rho_{\lambda_{2}}\)
  are quasi-compatible at~\(x\).
  (If \(\rho_{\Lambda}\) is a strongly quasi-compatible system,
  then this is automatic.
  If \(\rho_{\Lambda}\) is merely a quasi-compatible system,
  then this is true for \(x \in U^{\cl}\),
  for some non-empty open subset \(U \subset X\).)

  There exists an integer \(n \ge 1\)
  such that \(P(t) = P_{x,\rho_{\lambda},n}(t)\) does not depend on
  \(\lambda \in \Lambda\)
  (since we assumed that \(\Lambda\) is a finite set).
  We may then compute
  \[
   P_{x,\rho_{\lambda'},n}(t) =
   \prod_{\lambda|\lambda'}
   \Nm^{E_{\lambda}}_{E'_{\lambda'}} P_{x,\rho_{\lambda},n}(t)
   = \Nm^{E}_{E'} P(t).
  \]
  We conclude that
  \(P_{x,\rho_{\lambda'},n}(t)\) is a polynomial in~\(E'[t]\)
  that does not depend on \(\lambda' \in \Lambda'\).
 \end{proof}
\end{lemma}

\paragraph{} \label{tildeE} 
A counterpart to the previous lemma is as follows.
Let \(K\) be a finitely generated field.
Let \(E\) be a number field, and
let \(\Lambda\) be a set of finite places of~\(E\).
Let \(\rho_{\Lambda}\) be a system of Galois representations over~\(K\).
Let \(E \subset \tilde{E}\) be a finite extension,
and let \(\tilde{\Lambda}\) be
the set of finite places~\(\tilde{\lambda}\) of~\(\tilde{E}\)
that lie above places \(\lambda \in \Lambda\).

Let \(\lambda \in \Lambda\) be a finite place of~\(E\).
Write \(\tilde{E}_{\lambda}\) for \(\tilde{E} \otimes_{E} E_{\lambda}\)
and recall that \(\tilde{E}_{\lambda} =
 \prod_{\tilde{\lambda}|\lambda} \tilde{E}_{\tilde{\lambda}}\).
Consider the representation \(\tilde{\rho}_{\lambda} =
 \rho_{\lambda} \otimes_{E_{\lambda}} \tilde{E}_{\lambda}\),
and observe that it naturally decomposes as
\(\tilde{\rho}_{\lambda} =
 \bigoplus_{\tilde{\lambda}|\lambda} \tilde{\rho}_{\tilde{\lambda}}\),
where \(\tilde{\rho}_{\tilde{\lambda}} =
 \rho_{\lambda} \otimes_{E_{\lambda}} \tilde{E}_{\tilde{\lambda}}\).
We assemble these Galois representations~\(\tilde{\rho}_{\tilde{\lambda}}\)
in a system of Galois representations that we denote with
\(\tilde{\rho}_{\tilde{\Lambda}}\) or \(\rho_{\Lambda} \otimes_{E} \tilde{E}\).

\begin{lemma} \label{tildeElem} 
 Let \(K\), \(E \subset \tilde{E}\), \(\Lambda\), \(\tilde{\Lambda}\),
 \(\rho_{\Lambda}\), and~\(\tilde{\rho}_{\tilde{\Lambda}}\)
 be as in \cref{tildeE}.
 If \(\rho_{\Lambda}\) is
 a (strongly) quasi-compatible system of Galois representations,
 then \(\tilde{\rho}_{\tilde{\Lambda}}\) is
 a (strongly) quasi-compatible system of Galois representations.
 \begin{proof}
  Let \(X\) be a model of~\(K\) and let \(x \in X^{\cl}\) be a closed point.
  Let \(\tilde{\lambda} \in \tilde{\Lambda}\) be a place
  that lies above \(\lambda \in \Lambda\), and let \(n \ge 1\) be an integer.
  Then \(P_{x,\rho_{\lambda},n} = P_{x,\tilde{\rho}_{\tilde{\lambda}},n}\).
 \end{proof}
\end{lemma}

\begin{lemma} \label{plethysm} 
 Let \(K\) be a finitely generated field.
 Let \(E\) be a number field;
 and let \(\Lambda\) be a set of finite places of~\(E\).
 Let \(\rho_{\Lambda}\) and \(\rho'_{\Lambda}\)
 be two systems of Galois representations over~\(K\).
 Then one may naturally form the following systems of Galois representations:
 \begin{enumerate}[label=(\alph*),ref=(\alph*)]
  \item \label{pleth_first}
   the dual: \(\dual{\rho}_{\Lambda} =
   (E, \Lambda, (\dual{\rho})_{\lambda \in \Lambda})\);
  \item the direct sum:
   \(\rho_{\Lambda} \oplus \rho'_{\Lambda} =
   (E,\Lambda, (\rho_{\lambda} \oplus \rho'_{\lambda})_{\lambda \in \Lambda})\);
  \item the tensor product:
   \(\rho_{\Lambda} \otimes \rho'_{\Lambda} =
   (E,\Lambda,(\rho_{\lambda} \otimes \rho'_{\lambda})_{\lambda \in \Lambda})\);
  \item \label{pleth_last}
   the internal Hom:
   \(\iHom(\rho_{\Lambda}, \rho'_{\Lambda}) =
   (E, \Lambda,(\iHom(\rho_{\lambda},\rho'_{\lambda}))_{\lambda \in \Lambda})\).
 \end{enumerate}
 If \(\rho_{\Lambda}\) and \(\rho'_{\Lambda}\) are
 systems of Galois representations over~\(K\) that are both quasi-compatible,
 then the constructions~\labelcref{pleth_first} through~\labelcref{pleth_last}
 form a quasi-compatible system of Galois representations.
 \begin{proof}
  This \namecref{plethysm} is an immediate consequence
  of the following \namecref{eigenSchur}.
 \end{proof}
\end{lemma}

\begin{lemma} \label{eigenSchur} 
 Let \(V\) and~\(V'\) be finite-dimensional vector spaces over a field~\(K\).
 Let \(g\) and~\(g'\) be endomorphisms of~\(V\) and~\(V'\) respectively.
 \begin{enumerate}
  \item The coefficients of the characteristic polynomial
   \(\cp{g \oplus g' | V \oplus V'}\)
   are integral polynomial expressions in the coefficients of
   \(\cp{g | V}\) and \(\cp{g' | V'}\).
  \item The coefficients of \(\cp{g \otimes g' | V \otimes V'}\)
   are integral polynomial expressions in the coefficients of
   \(\cp{g | V}\) and \(\cp{g' | V'}\).
 \end{enumerate}
 \begin{proof}
  Write \(f\) for \(\cp{g | V}\) and \(f'\) for \(\cp{g'|V'}\).
  For point~1, note that
  \(\cp{g \oplus g' | V \oplus V} = f \cdot f'\).
  For point~2, put \(f = \prod_{i = 1}^{n}(x - \alpha_{i})\)
  and \(f' = \prod_{j = 1}^{n'}(x - \alpha'_{j})\) in \(\bar{K}[x]\),
  and note that
  \begin{align*}
   \cp{g \otimes g' | V \otimes V'}
   &= \prod_{i,j} (x - \alpha_{i}\alpha'_{j}) \\
   &= \prod_{i = 1}^{n} \alpha_{i}^{n'}
   \prod_{j = 1}^{n'} (x/\alpha_{i} - \alpha'_{j}) \\
   &= \res_{y}(f(y), f'(x/y) \cdot y^{n'}),
  \end{align*}
  where \(\res_{y}(\_,\_)\) denotes the resultant of the polynomials in~\(y\).
 \end{proof}
\end{lemma}

\section{Examples of quasi-compatible systems} 
\label{CSReg}

\readme 
In this section we show that
abelian varieties and abelian \cm~motives
give rise to strongly quasi-compatible systems of Galois representations
(in respectively \cref{abvarCSR} and \cref{cmMotCSR}).
These results are known over number fields.
We recall their proofs and generalise the results to finitely generated fields.

\paragraph{} 
Let \(K \subset \CC\) be a finitely generated field.
Let \(M\) be a motive over~\(K\).
Let \(E \subset \End(M)\) be a number field.
Let \(\Lambda\) be the set of finite places of~\(E\).
Let \(\ell\) be a prime number.
Then \(\Hl(M)\) is a module over
\(E_{\ell} = E \otimes \QQl \cong \prod_{\lambda|\ell} E_{\lambda}\).
Correspondingly, the Galois representation \(\Hl(M)\) decomposes as
\(\Hl(M) \cong \bigoplus_{\lambda|\ell} \Hlambda(M)\),
with \(\Hlambda(M) = \Hl(M) \otimes_{E_{\ell}} E_{\lambda}\).
The \(\lambda\)-adic representations \(\Hlambda(M)\),
with \(\lambda \in \Lambda\),
form a system of Galois representations
that we denote with \(\HH_{\Lambda}(M)\).
It is expected that \(\HH_{\Lambda}(M)\)
is a quasi-compatible system of Galois representations,
and even a compatible system in the sense of Serre.
(Indeed, this assertion is implied by the Tate conjecture.)

The following \namecref{abvarCSR} is a slightly weaker version of
a result proven by Shimura in~\S11.10.1 of~\cite{Sh67}.
We present the proof by Shimura in modern notation, and with a bit more detail.
The proof is given in \cref{abvarCSRprf}, and relies on \cref{lemShimura},
which is proposition~11.9 of~\cite{Sh67}.
For similar discussions, see \cite{Ch92}, \S\textsc{ii} of~\cite{Ribet},
\cite{No09}, and~\cite{No13}.

\begin{theorem}[\S11.10.1 of~\cite{Sh67}] \label{abvarCSR} 
 Let \(A\) be an abelian variety over a finitely generated field~\(K\);
 and let \(E \subset \End^{0}(A)\) be a number field.
 Let \(\Lambda\) be the set of finite places of~\(E\)
 whose residue characteristic is different from~\(\chrc(K)\).
 Then \(\HH^{1}_{\Lambda}(A)\) is
 a strongly quasi-compatible system of Galois representations.
 \begin{proof}
  See \cref{abvarCSRprf}.
 \end{proof}
\end{theorem}

\begin{proposition}[11.9 of~\cite{Sh67}] \label{lemShimura} 
 Let \(E\) be a number field.
 Let \(\primes\) be a set of prime numbers.
 Let \(\Lambda\) be the set of finite places of~\(E\)
 that lie above a prime number in~\(\primes\).
 For every prime number~\(\ell \in \primes\),
 let \(\Hl\) be a finitely generated \(E_{\ell}\)-module.
 (Recall that
 \(E_{\ell} = E \otimes \QQl \cong \prod_{\lambda|\ell} E_{\lambda}\).)
 Write \(\Hlambda\) for \(\Hl \otimes_{E_{\ell}} E_{\lambda}\),
 so that \(\Hl \cong \bigoplus_{\lambda|\ell} \Hlambda\).

 Let \(R\) be a finite-dimensional commutative semisimple \(E\)-algebra;
 and suppose that, for every prime number~\(\ell \in \primes\),
 we are given \(E\)-algebra homomorphisms \(R \to \End_{E_{\ell}}(\Hl)\).
 Assume that for every \(r \in R\)
 the characteristic polynomial \(\cp[\QQl]{r|\Hl}\)
 has coefficients in~\(\QQ\) and is independent of~\(\ell \in \primes\).
 Under these assumptions,
 for every \(r \in R\) the characteristic polynomial
 \(\cp[E_{\lambda}]{r |\Hlambda}\)
 has coefficients in~\(E\) and is independent of~\(\lambda \in \Lambda\).
 \begin{proof}
  The assumptions on \(R\) imply that~\(R\)
  is a finite product of finite field extensions~\(K_{i}/E\).
  Let \(\epsilon_{i}\) be the idempotent of~\(R\)
  that is~\(1\) on~\(K_{i}\) and~\(0\) elsewhere.
  For \(r \in R\), observe that
  \[
   \cp[\QQl]{r |\Hl} =
   \prod_{i} \cp[\QQl]{\epsilon_{i}r | \epsilon_{i}\Hl}, \qquad
   \cp[E_{\lambda}]{r|\Hlambda} =
   \prod_{i} \cp[E_{\lambda}]{\epsilon_{i}r | \epsilon_{i}\Hlambda}.
  \]
  We conclude that we only need to prove the lemma for \(R = K_{i}\),
  and \(\Hl = \epsilon_{i}\Hl\),
  \textit{i.e.,} that we can reduce to the case where \(R\) is a field.

  Suppose \(R\) is a finite field extension of~\(E\),
  and choose an element \(\pi \in R\) that generates \(R\) as a field.
  Let \(f^{\pi}_{\QQ}\) be the minimum polynomial of~\(\pi\) over~\(\QQ\).
  Observe that \(\cp[\QQl]{\pi |\Hl}\) is
  a divisor of a power of~\(f^{\pi}_{\QQ}\) in~\(\QQl[t]\).
  Since both are elements of~\(\QQ[t]\)
  and~\(f^{\pi}_{\QQ}\) is irreducible,
  we conclude that \(\cp[\QQl]{\pi |\Hl}\) is equal
  to~\((f^{\pi}_{\QQ})^{d}\), for some positive integer~\(d\).
  Since~\(\pi\) is semisimple,
  it follows that \(\Hl \cong \QQl[\pi]^{d}\) as \(\QQl[\pi]\)-modules.
  Let \(\HH\) be the \(R\)-vector space~\(R^{d}\).
  By construction \(\Hl \cong \HH \otimes_{\QQ} \QQl\)
  as \((R \otimes_{\QQ} \QQl)\)-modules.
  Because \(R \otimes_{\QQ} \QQl \cong R \otimes_{E} E \otimes_{\QQ} \QQl\),
  this implies that \(\Hlambda \cong \HH \otimes_{E} E_{\lambda}\)
  as \((R \otimes_{E} E_{\lambda})\)-modules.
  For all \(r \in R\), we have
  \(\cp[E_{\lambda}]{r|\Hlambda} = \cp[E]{r|\HH}\),
  and therefore \(\cp[E_{\lambda}]{r|\Hlambda}\)
  has coefficients in~\(E\) and is independent of~\(\lambda \in \Lambda\).
 \end{proof}
\end{proposition}

\begin{lemma} \label{abvarCSRfinfld} 
 Let \(A\) be an abelian variety
 over a finite field~\(\kappa\) of characteristic~\(p\).
 Let \(E\) be a number field inside \(\End^{0}(A)\).
 Let \(\Lambda\) be the set of finite places of~\(E\)
 whose residue characteristic is different from~\(p\).
 Then \(P_{x,\rho_{\lambda},1}(t)\) has coefficients in~\(E\)
 and is independent of~\(\lambda \in \Lambda\).
 In particular, \(\HH^{1}_{\Lambda}(A)\) is
 a quasi-compatible system of Galois representations.
 \begin{proof}
  Note that \(\Spec(\kappa)\) is the only model of~\(\kappa\).
  Let \(x\) denote the single point of~\(\Spec(\kappa)\).
  Let \(E[F_{x}]\) be the subalgebra of \(\End^{0}(A)\)
  generated by \(E\) and \(F_{\bar{\kappa}/\kappa}\).
  Note that \(E[F_{x}]\) may naturally be viewed as the subalgebra of
  \(\End(\Hl^{1}(A))\) generated by \(E\) and \(F_{x,\rho_{\ell}}\).
  This algebra is semisimple by work of~Weil.
  For every \(r \in E[F_{x}]\)
  the characteristic polynomial \(\cp{r | \Hl^{1}(A)}\)
  has coefficients in~\(\QQ\),
  and is independent of~\(\ell\), by theorem~2.2 of~\cite{KM74}.
  It follows from \cref{lemShimura} that
  \(P_{x,\rho_{\lambda},1}(t)\) has coefficients in~\(E\)
  and is independent of~\(\lambda \in \Lambda\).
 \end{proof}
\end{lemma}

\begin{corollary}[theorem~\textsc{ii}.2.1.1 of~\cite{Ribet}] \label{freeElmod} 
 Let \(A\) be an abelian variety over a finitely generated field~\(K\),
 and fix a prime number \(\ell \ne \chrc(K)\).
 Let \(E \subset \End^{0}(A)\) be a number field.
 Then \(\Hl^{1}(A)\) is a free \(E_{\ell}\)-module.
 \begin{proof}
  Let \(X\) be a model of~\(K\),
  and let \(x \in X^{\cl}\) be a closed point
  whose residue characteristic is different from~\(\ell\)
  and such that~\(A\) has good reduction at~\(x\).
  Specialise to~\(x\) and apply \cref{abvarCSRfinfld}.
 \end{proof}
\end{corollary}

\begin{lemma} \label{semiab} 
 Let \(K\) be a field.
 Let \(T \into G \stackrel{\alpha}{\longto} A\)
 be a semiabelian variety over~\(K\).
 Let \(E\) be a number field inside \(\End^{0}(G)\).
 Then \(E\) naturally maps to \(\End^{0}(A)\) and~\(\End^{0}(T)\).
 \begin{proof}
  Let \(f\) be an endomorphism of~\(G\).
  Consider the composition
  \(g \colon T \into G \stackrel{f}{\longto} G \onto A\).
  The image of \(g\) is affine,
  since it is a quotient of~\(T\),
  and it is projective, since it is a closed subgroup of~\(A\).
  It is also connected and reduced, and therefore factors via \(0 \in A(K)\).
  We conclude that \(f(T) \subset T\), which proves the result.
 \end{proof}
\end{lemma}

\begin{lemma} \label{NOS} 
 Let \(X\) be the spectrum of a discrete valuation ring.
 Let \(\eta\) (resp.~\(x\)) denote the generic (resp.\ special) point of~\(X\).
 Let \(A\) be a semistable abelian variety over~\(\eta\).
 Let \(E\) be a number field inside \(\End^{0}(A)\).
 Let \(\lambda\) be a finite place of~\(E\)
 such that the residue characteristics of~\(\lambda\) and~\(x\) are different.
 Then \(A\) has good reduction at~\(x\)
 if and only if \(\HH^{1}_{\lambda}(A)\) is unramified at~\(x\).
 \begin{proof}
  This is a slight generalisation of
  the criterion of N\'{e}ron--Ogg--Shafarevic,
  theorem~1 of~\cite{ST68}.
  It is clear that if \(A\) has good reduction at~\(x\),
  then \(\HH^{1}_{\lambda}(A)\) is unramified at~\(x\).
  We focus on the converse implication.
  Let \(\ell\) be the residue characteristic of~\(\lambda\).
  By theorem~1 of~\cite{ST68} it suffices to show that
  \(\Hl^{1}(A)\) is unramified at~\(x\).
  Let \(\Hl^{1}(A)^{I}\) denote the subspace of
  \(\Hl^{1}(A)\) that is invariant under inertia.
  Let \(G\) be the N\'{e}ron model of~\(A\) over~\(X\).
  Recall that \(\Hl^{1}(A)^{I} \cong \Hl^{1}(G_{x})\),
  by lemma~2 of~\cite{ST68}.
  It follows from the definition of the N\'{e}ron model
  that \(E\) embeds into \(\End(G) \otimes \QQ\).
  Hence \(E\) embeds into \(\End(G_{x}) \otimes \QQ\), and
  we claim that \(\Hl^{1}(A)^{I} \cong \Hl^{1}(G_{x})\)
  is a free \(E_{\ell}\)-module.
  (With \(E_{\ell}\) we mean
  \(E \otimes \QQl \cong \prod_{\lambda | \ell} E_{\lambda}\).)
  Before proving the claim, let us see why it is sufficient
  for proving the lemma.
  By \cref{freeElmod} we know that \(\Hl^{1}(A)\)
  is a free \(E_{\ell}\)-module.
  Thus \(\Hl^{1}(A)/\Hl^{1}(A)^{I}\) is a free \(E_{\ell}\)-module.
  We conclude that \(\HH^{1}_{\lambda}(A)\) is unramified at~\(x\),
  if and only if \(\Hl^{1}(A)\) is unramified at~\(x\).

  We will now prove the claim that \(\Hl^{1}(A)^{I} \cong \Hl^{1}(G_{x})\)
  is a free \(E_{\ell}\)-module.
  Since \(A\) is semistable, the special fibre \(G_{x}\)
  is a semiabelian variety \(T \into G_{x} \to B\).
  The semiabelian variety \(G_{x}\) is a special case of a \(1\)-motive,
  and thus we have a short exact sequence
  \[
   0 \to \Hl^{1}(B) \to \Hl^{1}(G_{x}) \to \Hl^{1}(T) \to 0.
  \]
  We also have \(\Hl^{1}(T) \cong \Hom(T,\Gm) \otimes \QQl(-1)\),
  see variante~10.1.10 of~\cite{HodgeIII}.
  By \cref{semiab}, the action of~\(E\) on~\(G_{x}\)
  gives an action of~\(E\) on both~\(T\) and~\(B\).
  Since \(\Hom(T,\Gm) \otimes \QQ\) is a free \(E\)-module,
  we know that \(\Hl^{1}(T)\) is a free \(E_{\ell}\)-module.
  By \cref{freeElmod} we also know that \(\Hl^{1}(B)\)
  is free as \(E_{\ell}\)-module.
  Therefore, \(\Hl^{1}(G_{x}) \cong \Hl^{1}(A)^{I}\)
  is free as \(E_{\ell}\)-module.
 \end{proof}
\end{lemma}

\begin{nproof}[of \cref{abvarCSR}] \label{abvarCSRprf} 
 Let \(X\) be a model of~\(K\); and let \(x \in X^{\cl}\) be a closed point.
 Let \(\Lambda^{(x)}\) be the set of places \(\lambda \in \Lambda\)
 that have a residue characteristic~\(\ell\) that is different
 from the residue characteristic of~\(x\).
 If there is a \(\lambda \in \Lambda^{(x)}\)
 such that \(\HH^{1}_{\lambda}(A)\) is unramified at~\(x\),
 then \(A\) has good reduction at~\(x\), by \cref{NOS}.
 Assume that \(A\) has good reduction at~\(x\).
 We denote this reduction with~\(A_{x}\).
 It follows from \cref{abvarCSRfinfld} that
 \(P_{x,\rho_{\lambda},1}(t)\) has coefficients in~\(E\)
 and is independent of~\(\lambda \in \Lambda^{(x)}\).
\end{nproof}

\begin{proposition} \label{halftwistCSR} 
 Let \(M\) be an abelian motive of weight~\(n\)
 over a finitely generated field \(K \subset \CC\).
 Let \(E \subset \End(M)\) be a \cm~field~\(E\)
 such that \(\dim_{E}(M) = 1\),
 and let \(\Lambda\) be the set of finite places of~\(E\).
 Then the system \(\HH_{\Lambda}(M)\) is a
 strongly quasi-compatible system of Galois representations.
 \begin{proof}
  Let \(m\) be the level of~\(M\),
  that is \(\max\{p-q \mid \HB(M)^{p,q} \ne 0\}\).
  We apply induction to~\(m\),
  and use half-twists as described in \cref{halftwists}.
  If \(m = 0\), then there is nothing to be done.

  Suppose that \(m \ge 1\).
  Let \(T \subset \Sigma(E)\) be the set of embeddings
  through which \(E\) acts on
  \(\bigoplus_{p \ge \lceil n/2\rceil}\HB(M)^{p,q}\).
  Since \(\dim_{E}(M) = 1\) we know that \(T \cap T^{\dagger} = \varnothing\).
  It follows from the discussion in \cref{halftwists} and \cref{halftwistsMot}
  that there exists a finitely generated extension \(L/K\),
  an abelian variety~\(A\) over~\(L\), and a motive~\(N\) over~\(L\) such that
  \(M_{L} \cong \iHom_{E}(\HH^{1}(A), N)\),
  and such that the level of~\(N\) is \(m-1\), and \(\dim_{E}(N) = 1\).
  By \cref{abvarCSR} we know that
  \(\HH^{1}_{\Lambda}(A)\) is a strongly quasi-compatible system,
  and by induction we may assume that \(\HH^{1}_{\Lambda}(N)\)
  is a strongly quasi-compatible system.
  It follows from \cref{plethysm} that
  \(\HH_{\Lambda}(M_{L}) \cong
  \iHom_{E}(\HH^{1}_{\Lambda}(A), \HH^{1}_{\Lambda}(N))\)
  is a quasi-compatible system of Galois representations over~\(L\),
  and we will now argue that it is even a strongly quasi-compatible system.

  Let \(X\) be a model of~\(K\), and let~\(x \in X\) be a closed point.
  Let \(\Lambda^{(x)}\) be the set of finite places of~\(E\)
  whose residue characteristic is
  different from the residue characteristic of~\(x\).
  Fix \(\lambda \in \Lambda^{(x)}\).
  We may assume that \(A\) is semistable over~\(L\)
  (possibly replacing~\(L\) with a finite field extension).
  Since \(A\) is a semistable \cm~abelian variety,
  we know that \(A\) has good reduction everywhere,
  and thus \(\HH^{1}_{\lambda}(A)\) is unramified at~\(x\).
  Hence \(\HH_{\lambda}(M_{L})\) is unramified at~\(x\)
  if and only if \(\HH^{1}_{\lambda}(N)\) is unramified at~\(x\).
  Finally, \cref{CSRbc} shows that \(\HH_{\Lambda}(M)\) is also
  a strongly quasi-compatible system of Galois representations over~\(K\).
 \end{proof}
\end{proposition}

\begin{theorem}[see also corollary~\textsc{i}.6.5.7 of~\cite{Sch88}] 
 \label{cmMotCSR}
 Let~\(M\) be an abelian \cm~motive
 over a finitely generated field~\(K \subset \CC\).
 Let \(E\) be a subfield of~\(\End(M)\),
 and let \(\Lambda\) be the set of finite places of~\(E\).
 Then the system \(\HH_{\Lambda}(M)\) is
 a strongly quasi-compatible system of Galois representations.
 \begin{proof}
  By \cref{CSRbc} we may replace~\(K\) by a finitely generated extension
  and thus we may and do assume that \(M\)
  decomposes into a sum of geometrically isotypical components
  \(M = M_{1} \oplus \ldots \oplus M_{r}\).
  Observe that \(E \subset \End(M_{i})\) for \(i = 1,\ldots,r\).
  By \cref{plethysm} we see that it suffices to show that
  \(\HH_{\Lambda}(M_{i})\) is a strongly quasi-compatible system
  for \(i = 1,\ldots,r\).
  Therefore we may assume that \(M \cong (M')^{\oplus k}\),
  where \(M'\) is a geometrically irreducible \cm-motive.
  If \(M'\) is a Tate motive,
  then the result is trivialy true.
  Hence, let us assume that \(E' = \End(M')\) is a \cm~field.
  Notice that \(\dim_{E'}(M') = 1\).
  By assumption \(E\) acts on \((M')^{\oplus k}\),
  and thus we get a specific embedding \(E \subset \Mat_{k}(E')\).
  We may find a field~\(\tilde{E} \subset \Mat_{k}(E')\)
  that contains \(E\),
  and such that \([\tilde{E} : E'] = k\).
  Then \(M = M' \otimes_{E'} \tilde{E}\).
  Let \(\tilde{\Lambda}\) be the set of finite places of~\(\tilde{E}\).
  (N.b., we now have inclusions \(E' \subset E \subset \tilde{E}\).)
  By \cref{halftwistCSR}, the system
  \(\HH_{\Lambda'}(M')\)
  is a strongly quasi-compatible (\(E'\)-rational) system
  of Galois representations,
  and by \cref{tildeElem} we find that
  \(\HH_{\tilde{\Lambda}}(M) = \HH_{\Lambda'}(M') \otimes_{E'} \tilde{E}\)
  is a strongly quasi-compatible (\(\tilde{E}\)-rational) system.
  We conclude that \(\HH_{\Lambda}(M)\)
  is a strongly quasi-compatible (\(E\)-rational) system
  of Galois representations
  by \cref{restrictElem}.
  \looseness=-1
 \end{proof}
\end{theorem}

\section{Deformations of abelian motives} \label{mainproof} 

\readme 
In this section we prove the main result of this article,
which is the following \namecref{abmotCSR}.

\begin{theorem} \label{abmotCSR} 
Let \(M\) be an abelian motive
over a finitely generated field~\(K \subset \CC\).
Let \(E\) be a subfield of~\(\End(M)\),
and let \(\Lambda\) be the set of finite places of~\(E\).
Then the system \(\HH_{\Lambda}(M)\) is
a quasi-compatible system of Galois representations.
\end{theorem}

\paragraph{} 
The proof of this theorem relies heavily on the fact
that an abelian motive can be placed as fibre in a
family of abelian motives over a Shimura variety of Hodge type.
\Cref{biglist} summarises this result.
Its proof uses the rather technical \cref{fundconstr}.
Once we have the family of motives in place,
the rest of the section is devoted to the proof of the main theorem.
The following picture aims to capture the intuition of the proof.
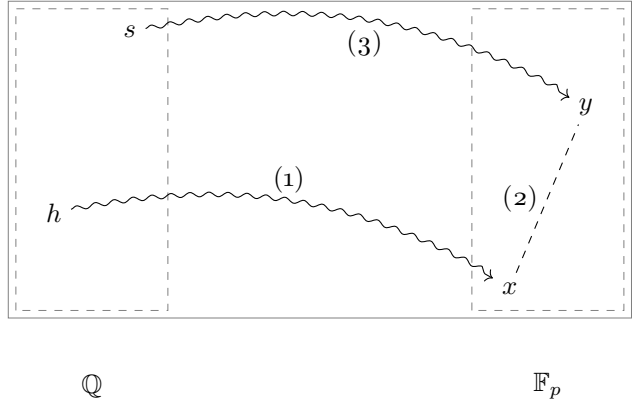
\[ 
 \begin{tikzpicture}
  \draw[help lines] (-.1,-.1) rectangle (8.1,4.1);
  \draw[help lines,dashed] (0,0) rectangle (2,4);
  \draw[help lines,dashed] (6,0) rectangle (8,4);
  \node at (1,-1) {$\QQ$};
  \node at (7,-1) {$\FFp$};

  \node (h) at (0.5,1.3) {$h$};
  \node (x) at (6.5,0.3) {$x$};
  \node (y) at (7.5,2.7) {$y$};
  \node (s) at (1.5,3.7) {$s$};

  \draw (h) edge[->,bend left=20,decorate,
  decoration={snake,amplitude=.3mm,segment length=2.5mm,post length=.5mm}]
  node[above] {\hbox{(1)}} (x);
  \draw (x) edge[dashed] node[left] {\hbox{(2)}} (y);
  \draw (s) edge[->,bend left=20,decorate,
  decoration={snake,amplitude=.3mm,segment length=2.5mm,post length=.5mm}]
  node[below] {\hbox{(3)}} (y);
 \end{tikzpicture}
\] 
The picture is a cartoon of an integral model of a Shimura variety,
and the motive~\(M\) fits into a family~\(\mathcal{M}\) over the generic fibre,
such that \(M \cong \mathcal{M}_{h}\).
We give a rough sketch of the strategy for the proof
that explains the three steps in the picture:
\begin{enumerate*}[label=(\arabic*)]
\item We have a system of Galois representations
 \(\HH_{\Lambda}(\mathcal{M}_{h})\)
 and we want to show that it is quasi-compatible at~\(x\);
\item we replace~\(x\) by an isogenous point~\(y\)
 (in the sense of Kisin~\cite{Kisin_modp}); and
\item we may assume that~\(y\) lifts to a special point~\(s\).
\end{enumerate*}
The upshot is that we have to show that the system
\(\HH_{\Lambda}(\mathcal{M}_{s})\) is quasi-compatible at~\(y\).
We will see that this follows from \cref{cmMotCSR}.

\begin{construction} \label{fundconstr} 
 Fix an integer \(g \in \ZZ_{\ge 0}\).
 Let \((G,X) \into (\GSp_{2g}, \Sds)\) be a morphism of Shimura data,
 and let \(h \in X\) be a morphism \(\DelS \to G_{\RR}\).
 In this paragraph we will construct
 an abelian scheme over an integral model 
 of the Shimura variety \(\Sh_{\Kmpt}(G,X)\),
 where \(\Kmpt\) is a certain compact open subgroup of~\(G(\Adele_{\fin})\).
 Along the way, we make two choices,
 labeled (\textit{i}) and~(\textit{ii}) so that we may refer to them later on.
 
 For each integer \(n \ge 3\), let \(\Kmpt'_{(n)}\) denote
 the congruence subgroup of \(\GSp_{2g}(\hat{\ZZ})\)
 consisting of elements congruent to~\(1\) modulo~\(n\).
 Write \(\Kmpt_{(n)}\) for \(\Kmpt'_{(n)} \cap G(\Adele_{\fin})\).
 This gives a morphism of Shimura varieties
 \[
  \Sh_{\Kmpt_{(n)}}(G,X) \to \Sh_{\Kmpt'_{(n)}}(\GSp_{2g},\Sds).
 \]
 By applying lemma~3.3 of~\cite{No96} with \(p = 6\)
 we can choose~\(n\) in such a way that it is coprime with~\(p\)
 and such that this morphism of Shimura varieties is a closed immersion.
 (In \cite{No96}, Noot assumes that \(p\) is prime,
 but he does not use this fact in his proof.)
 
 \begin{enumerate*}[label=(\textit{\roman*})]
  \item Fix such an integer~\(n\),
 \end{enumerate*}
 and write \(\Kmpt\) for~\(\Kmpt_{(n)}\).
 Since \(n > 3\), the subgroup \(\Kmpt'_{(n)}\) is neat,
 hence \(\Kmpt\) is neat, and therefore \(\Sh_{\Kmpt}(G,X)\) is smooth.
 As is common, we denote with \(\AV_{g,1,n}/\ZZ[1/n]\) the moduli space of
 principally polarised abelian varieties
 of dimension~\(g\) with a level-\(n\) structure.
 Recall that \(\AV_{g,1,n}\) is smooth over \(\ZZ[1/n]\).
 
 We have a closed immersion of Shimura varieties
 \[
  \Sh_{\Kmpt}(G,X) \into \AV_{g,1,n,\CC}.
 \]
 Let \(F' \subset \CC\) be the reflex field of~\((G,X)\).
 Let \(\mSh_{\Kmpt}(G,X)\) be the Zariski closure of \(\Sh_{\Kmpt}(G,X)\)
 in \(\AV_{g,1,n}\) over~\(\mathcal{O}_{F'}[1/n]\).
 There exists an integer multiple~\(N_{0}\) of~\(n\)
 such that \(\mSh_{\Kmpt}(G,X)_{\mathcal{O}_{F'}[1/N_{0}]}\) is smooth.
 
 For a prime number~\(p\),
 let \(\Kmpt_{p}\) be \(\Kmpt \cap G(\QQp)\),
 and let \(\Kmpt^{p}\) be \(\Kmpt \cap G(\Adele^{p}_{\fin})\).
 The set of prime numbers for which \(\Kmpt \ne \Kmpt_{p}\Kmpt^{p}\) is finite.
 Write~\(N_{1}\) for the product of those prime numbers.
 Let \(p\) be a prime number that does not divide~\(N_{1}\).
 The group~\(\Kmpt_{p}\) is called \emph{hyperspecial}
 if there is a reductive model \(\mathcal{G}/\ZZ_{p}\) of~\(G/\QQ\)
 such that \(\Kmpt_{p} = \mathcal{G}(\ZZ_{p})\).
 The set of prime numbers for which \(\Kmpt_{p}\) is not hyperspecial is finite.
 Write \(N_{2}\) for the product of those prime numbers.
 Let \(N\) be the integer \(N_{0} \cdot N_{1} \cdot N_{2}\).
 (We need to assure that \(\Kmpt_{p}\) is hyperspecial, for \(p \nmid N\),
 in order to apply results by Kisin~\cite{Kisin_modp} in
 \cref{Kisin_isog} and \cref{isoglift}.)

 The point~\(h \in X\) is a complex point of~\(\mSh_{\Kmpt}(G,X)\).
 After replacing~\(F'\) by a finite extension \(F \subset \CC\)
 we may assume that the generic fibre of the irreducible component
 \(\mSh \subset \mSh_{\Kmpt}(G,X)_{\mathcal{O}_{F}[1/N]}\)
 that contains the point~\(h\) is geometrically irreducible.

 \begin{enumerate*}[label=(\textit{\roman*}),resume]
  \item Choose such a field \(F \subset \CC\).
 \end{enumerate*}
 In the following paragraphs
 we will consider the closed immersion of Shimura varieties
 \(\mSh \into \AV_{g,1,n}\)
 as a morphism of schemes over~\(\mathcal{O}_{F}[1/N]\).
\end{construction}

\begin{lemma} \label{biglist} 
 Let \(M\) be an abelian motive
 over a finitely generated field~\(K \subset \CC\).
 There exist
 \begin{itemize}
  \item finitely generated fields \(F \subset L \subset \CC\),
   with \(K \subset L\);
  \item a smooth irreducible component~\(\mSh\)
   of an integral model of a Shimura variety over~\(F\),
   such that the generic fibre~\(\mSh_{F}\) is geometrically irreducible;
  \item an abelian scheme \(f \colon \mathcal{A} \to \mSh\);
  \item an idempotent motivated cycle~\(\gamma\) in
   \(\End((\mathbb{T}^{a,b}\mathrm{R}^{1}f_{\CC,*} \QQ)(m))\),
   for certain \(a, b, m \in \ZZ\);
  \item a family of abelian motives~\(\mathcal{M}/\mSh_{L}\),
   such that \(\mathcal{M}/\mSh_{\CC} \cong \Im(\gamma)\);
  \item an isomorphism \(M_{L} \cong \mathcal{M}_{h}\),
   for some point \(h \in \mSh(L)\).
 \end{itemize}
 \begin{proof}
  Since~\(M\) is an abelian motive,
  there exists a principally polarised complex abelian variety~\(A\)
  such that \(M_{\CC} \in \Tangen{A}\).
  Write \(V\) for \(\HB(M)\).
  Observe that \(\GB(V)\) is naturally a quotient of~\(\GB(A)\).
  Write \(G\) for~\(\GB(A)\), and let~\(h \colon \DelS \to G_{\RR}\)
  be the map that defines the Hodge structure on~\(\HB(A)\).
  Let \(X\) be the \(G(\RR)\)-orbit of~\(h\) in \(\Hom(\DelS, G_{\RR})\).
  Let \(g\) be \(\dim(A)\).
  The pair \((G,X)\) is a Shimura datum,
  and by construction we get a morphism of Shimura data
  \((G,X) \into (\GSp_{2g}, \Sds)\).
  Now run \cref{fundconstr}, choosing
  \begin{enumerate*}[label=(\textit{\roman*})]
   \item an integer~\(n\);
   \item a number field~\(F \subset \CC\);
  \end{enumerate*}
  and producing a closed immersion of Shimura varieties
  \(\mSh \into \AV_{g,1,n}\) over~\(\mathcal{O}_{F}[1/N]\).
  
  It follows from \cref{fundconstr},
  that the Hodge structure~\(V\) gives rise to
  a variation of Hodge structure~\(\mathcal{V}\)
  on \(\mSh_{\CC}\) such that the fibre of~\(\mathcal{V}\) above~\(h\) is~\(V\),
  and such that \(h\) is a Hodge generic point of \(\mSh_{\CC}\)
  with respect to the variation~\(\mathcal{V}\).
  The embedding \(\mSh \into \AV_{g,1,n}\)
  gives a natural abelian scheme \(f \colon \mathcal{A} \to \mSh\).
  The point~\(h\) is also a Hodge generic point with respect to~\(f\).
  Observe that \(A = \mathcal{A}_{h}\).

  Recall that \(V \in \Tangen{\HB^{1}(\mathcal{A}_{h})}\),
  which means that there exist integers~\(a\), \(b\), and \(m\),
  and some projector~\(\gamma_{h}\) on
  \(\mathbb{T}^{a,b}\HB^{1}(\mathcal{A}_{h})(m)\)
  whose image is isomorphic to~\(V\).
  Since \(h\) is a Hodge generic point of~\(\mSh_{\CC}\),
  the projector~\(\gamma_{h}\) spreads out to a projector \(\gamma\) on
  \((\mathbb{T}^{a,b}\mathrm{R}^{1}f_{\CC,*}\QQ)(m)\),
  and \(\mathcal{V}_{\mSh_{\CC}} \cong \Im(\gamma)\).

  By \cref{hdgmot}, the projector \(\gamma\) is motivated,
  and thus we obtain a family of abelian motives \(\mathcal{M}/\mSh_{\CC}\)
  whose Betti realisation is \(\mathcal{V}_{\mSh_{\CC}}\).
  In particular \(\mathcal{M}_{h} \cong M_{\CC}\).
  Finally, the point~\(h\), the projector~\(\gamma\),
  and the family of motives~\(\mathcal{M}\)
  are all defined over a finitely generated subfield \(L \subset \CC\)
  that contains \(F\) and~\(K\).
 \end{proof}
\end{lemma}

\paragraph{} \label{startCSRproof} 
We will now start the proof of \cref{abmotCSR}.
We retain the assumptions and notation of \cref{fundconstr} and \cref{biglist}.
Write \(S\) for \(\mSh_{L}\).
Let \(\mathcal{V}'\) be the variation of Hodge structure
\(\mathrm{R}^{1}f_{\CC,*}\QQ\) over~\(S(\CC)\),
and write~\(\mathcal{V}\) for the image of~\(\gamma\)
in~\((\mathbb{T}^{a,b}\mathcal{V}')(m)\);
it is a variation of Hodge structure that is the
Betti realisation of~\(\mathcal{M}/S(\CC)\).
Because~\(h\) is a Hodge generic point,
the field~\(E\) is a subfield of~\(\End(\mathcal{V})\).
Let \((e_{i})_{i}\) be a basis of~\(E\) as \(\QQ\)-vector space.

Let \(\ell\) be a prime number.
Let \(\mathcal{V}'_{\ell}\) be the lisse \(\ell\)-adic sheaf
\(\mathrm{R}^{1}f_{*}\QQl\) over~\(\mSh\).
By \cref{hdgmot},
the projector~\(\gamma\) on~\((\mathbb{T}^{a,b}\mathcal{V}')(m)\)
induces a projector
on~\((\mathbb{T}^{a,b}\mathcal{V}'_{\ell,S})(m)\) over~\(S\)
that spreads out to a projector~\(\gamma_{\ell}\)
on~\((\mathbb{T}^{a,b}\mathcal{V}'_{\ell})(m)\)
over the entirety of~\(\mSh\).
Let \(\mathcal{V}_{\ell}\) denote the image of~\(\gamma_{\ell}\).
Note that \(\mathcal{V}_{\ell,S}\)
is the \(\ell\)-adic realisation of~\(\mathcal{M}/S\).

By \cref{hdgmot} we see that~\(E_{\ell} = E \otimes \QQl\) is a subalgebra of
\(\End(\mathcal{V}_{\ell,S})\).
Since \(S\) is the generic fibre of~\(\mSh\),
we see that \(E_{\ell} \subset \End(\mathcal{V}_{\ell})\).
This has two implications, namely
\begin{enumerate*}[label=(\textit{\roman*})]
 \item we obtain classes \(e_{i,\ell} \in \End(\mathcal{V}_{\ell})\)
  that form a \(\QQl\)-basis for~\(E_{\ell}\); and
 \item because \(E_{\ell} = E \otimes \QQl
   \cong \prod_{\lambda | \ell} E_{\lambda}\),
  the lisse \(\ell\)-adic sheaf~\(\mathcal{V}_{\ell}\)
  decomposes as a sum \(\bigoplus_{\lambda | \ell} \mathcal{V}_{\lambda}\)
  of lisse \(\lambda\)-adic sheaves, where
  \(\mathcal{V}_{\lambda} = \mathcal{V}_{\ell} \otimes_{E_{\ell}} E_{\lambda}\).
\end{enumerate*}

\paragraph{} \label{Kisin_isog} 
Let \(p\) be a prime number that does not divide~\(N\),
so that \(\Kmpt\) decomposes as \(\Kmpt_{p}\Kmpt^{p}\),
and \(\Kmpt_{p}\) is hyperspecial.
Let \(\FFq/\FFp\) be a finite field.
Let \(x \in \mathscr{S}(\FFq)\) be a point.
Kisin defines the \emph{isogeny class of~\(x\)}
in \S1.4.14 of~\cite{Kisin_modp}.
It is a subset of \(\mathscr{S}(\FFqbar)\).

Let~\(y\) be a point in~\(\mathscr{S}(\FFqbar)\)
that is isogenous to~\(x\).
Proposition~1.4.15 of~\cite{Kisin_modp}
implies that there is
an isomorphism of Galois representations
\(\mathcal{V}'_{\ell,x} \cong \mathcal{V}'_{\ell,y}\)
such that
\(\gamma_{\ell,x} \in \End((\mathbb{T}^{a,b}\mathcal{V}'_{\ell,x})(m))\)
is mapped to
\(\gamma_{\ell,y} \in \End((\mathbb{T}^{a,b}\mathcal{V}'_{\ell,y})(m))\),
and such that~\(e_{i,\ell,x}\) is mapped to~\(e_{i,\ell,y}\).
This implies that \(\mathcal{V}_{\ell,x} \cong \mathcal{V}_{\ell,y}\)
as \(E_{\ell}[\Gal(\FFqbar/\FFq)]\)-modules.
We conclude that \(\mathcal{V}_{\lambda,x} \cong \mathcal{V}_{\lambda,y}\)
as \(\lambda\)-adic Galois representations.

\paragraph{} \label{isoglift} 
We need one more key result by Kisin~\cite{Kisin_modp}.
Theorem~2.2.3 of~\cite{Kisin_modp} states that for every point
\(x \in \mathscr{S}(\FFqbar)\),
there is a point \(y \in \mathscr{S}(\FFqbar)\)
that is isogenous to~\(x\)
and such that~\(y\) is the reduction of a special point in~\(S\).

\paragraph{} \label{CSRproof} 
We are now set for the attack on \cref{abmotCSR}.
Let \(\lambda_{1}\) and~\(\lambda_{2}\) be two finite places of~\(E\).
Let \(\ell_{1}\) and~\(\ell_{2}\) be the residue characteristics
of~\(\lambda_{1}\) respectively~\(\lambda_{2}\).
Let \(X\) be the Zariski closure of~\(h\) in~\(\mSh\).
Note that \(X\) is a model for the residue field of~\(h\).
Let \(U \subset X\) be the Zariski open locus of points \(x \in X\)
such that the residue characteristic~\(p\) of~\(x\)
does not divide \(N \cdot \ell_{1} \cdot \ell_{2}\).
To prove \cref{abmotCSR}, it suffices to show that
\(\HH_{\lambda_{1}}(M)\) and~\(\HH_{\lambda_{2}}(M)\)
are quasi-compatible at all points~\(x \in U^{\cl}\).
Fix a point \(x \in U^{\cl}\).
Observe that by construction the representations
\(\HH_{\lambda_{1}}(M)\) and~\(\HH_{\lambda_{2}}(M)\)
are unramified at~\(x\).
Let \(\FFq\) be the residue field of~\(x\).
We want to show that
\(\mathcal{V}_{\lambda_{1},x}\) and~\(\mathcal{V}_{\lambda_{2},x}\)
are quasi-compatible.
This means that we have to show that the characteristic polynomials
of the Frobenius endomorphisms of
\(\mathcal{V}_{\lambda_{1},x}\) and~\(\mathcal{V}_{\lambda_{2},x}\)
are equal, possibly after replacing the Frobenius endomorphism by some power.
Equivalently, we may pass to a finite extension of~\(\FFq\).
This is what we will now do.

As mentioned in~\cref{isoglift}, theorem~2.2.3 of~\cite{Kisin_modp}
shows that there exists a point~\(y \in \mSh(\FFqbar)\)
such that \(y\) is isogenous to~\(x\)
and such that \(y\) is the reduction of a special point~\(s \in S\).
The point~\(y\) is defined over a finite extension of~\(\FFq\).
As explained in the previous paragraph,
we may replace~\(\FFq\) with a finite extension.
Thus we may and do assume that \(y\) is \(\FFq\)-rational.

We want to prove that
\(\mathcal{V}_{\lambda_{1},x}\) and~\(\mathcal{V}_{\lambda_{2},x}\)
are quasi-compatible.
By our remarks in \cref{Kisin_isog} we may as well show that
\(\mathcal{V}_{\lambda_{1},y}\) and~\(\mathcal{V}_{\lambda_{2},y}\)
are quasi-compatible.
In other words, we may show that
\(\HH_{\lambda_{1}}(\mathcal{M}_{s})\)
and~\(\HH_{\lambda_{2}}(\mathcal{M}_{s})\)
are quasi-compatible at~\(y\).
Recall that~\(s\) is a special point in~\(S\).
Therefore \(\mathcal{M}_{s}\) is an abelian \cm~motive,
and we conclude by \cref{cmMotCSR} that
\(\HH_{\lambda_{1}}(\mathcal{M}_{s})\)
and~\(\HH_{\lambda_{2}}(\mathcal{M}_{s})\)
are quasi-compatible at~\(y\).
This completes the proof of \cref{abmotCSR}.

\begin{remark} \label{compLask} 
 Laskar~\cite{Lask14} has obtained similar results.
 Let \(M\) be an abelian motive over a number field~\(K\).
 Let \(E \subset \End(M)\) be a number field,
 and let \(\Lambda\) be the set of finite places of~\(E\).
 Laskar needs the following condition:
 Assume that \(\GB(M)^{\ad}\)
 does not have a simple factor whose Dynkin diagram has type~\(D_{k}\)
 (or, more precisely, type~\(D_{k}^{\HQ}\) in the sense of
 table~1.3.8 of~\cite{Del_ShimVar}).
 Then theorem~1.1 of~\cite{Lask14} implies that the system~\(\HH_{\Lambda}(M)\)
 is a compatible system in the sense of Serre
 (that is, one may take \(n = 1\) in \cref{comprep})
 after replacing \(K\) by a finite extension.
 If \(\GB(M)^{\ad}\) does have a simple factor
 whose Dynkin diagram has type~\(D_{k}^{\HQ}\),
 then Laskar also obtains results,
 but I do not see how to translate them into our terminology.
 See~\cite{Lask14} for more details.
\end{remark}

\section{Properties of quasi-compatible systems} \label{props} 

\readme 
We derive some useful properties of quasi-compatible systems:
\begin{enumerate}
 \item (\Cref{FrobTor})
  We recover the notion of a Frobenius torus,
  just as in the classical concept of a compatible system
  in the sense of Serre.
 \item (\Cref{CSRisom})
  Let \(K\) be a finitely generated field.
  Let \(E\) be a number field, and let \(\Lambda\)
  be a set of finite places of~\(E\).
  Let \(\rho_{\Lambda}\) and \(\rho'_{\Lambda}\)
  be two \(E\)-rational quasi-compatible system
  of semisimple Galois representations of~\(K\).
  If there is a place \(\lambda \in \Lambda\)
  such that \(\rho_{\lambda} \cong \rho'_{\lambda}\) as
  \(\lambda\)-adic Galois representations,
  then \(\rho_{\Lambda} \cong \rho'_{\Lambda}\)
  as \(E\)-rational systems of Galois representations of~\(K\).
 \item (\Cref{recoverEnd})
  We show that under reasonable conditions,
  we can recover the field~\(E\) as subring of
  \(\End_{\Gal(\bar{K}/K),E_{\lambda}}(\rho_{\lambda})\)
  for some \(\lambda \in \Lambda\).
 \item (\Cref{rkindeplambda})
  Let \(\rho_{\Lambda}\) be a quasi-compatible system of
  semisimple Galois representations.
  We prove that the rank of \(\GG_{\lambda}(\rho_{\lambda})\)
  is independant of~\(\lambda\).
\end{enumerate}

\begin{definition}[see also \S3 of~\cite{Ch92}] \label{FrobTor} 
 Let \(K\) be a finitely generated field,
 let \(X\) be a model of~\(K\), and let \(x \in X^{\cl}\) be a closed point.
 Let \(E\) be a number field, and let \(\lambda\) be a finite place of~\(E\).
 Let \(\rho\) be a semisimple \(\lambda\)-adic Galois representation of~\(K\).
 Assume that \(\rho\) is unramified at~\(x\).
 The algebraic subgroup \(H_{n} \subset \GG_{\lambda}(\rho)\)
 generated by \(F_{x,\rho}^{n}\) is well-defined up to conjugation.
 Note that \(H_{n}\) is a finite-index subgroup of \(H_{1}\),
 and therefore the identity component of~\(H_{n}\)
 does not depend on~\(n\).
 We denote this identity component with \(T_{x}(\rho)\).

 If there is an integer \(n > 0\) such that \(F_{x,\rho}^{n}\) 
 is semisimple, then we call \(T_{x}(\rho)\)
 the \emph{Frobenius torus} at~\(x\).
 In this case the algebraic group \(T_{x}(\rho)\) is indeed an algebraic torus,
 which means that \(T_{x}(\rho)_{\bar{E}_{\lambda}} \cong \Gm^{k}\),
 for some \(k \ge 0\).
\end{definition}

\begin{remark} \label{FrobTorE} 
 Let \(K\) be a finitely generated field,
 let \(X\) be a model of~\(K\), and let \(x \in X^{\cl}\) be a closed point.
 Let \(E\) be a number field, and let \(\lambda\) be a finite place of~\(E\).
 Let \(\rho\) be a semisimple \(\lambda\)-adic Galois representation of~\(K\).
 Assume that there is an integer \(n > 0\)
 such that \(F_{x,\rho}^{n}\) is semisimple.

 Now assume that \(\rho\) is \(E\)-rational.
 Fix an integer \(n > 0\) such that \(F_{x,\rho}^{n}\) is semisimple
 and generates the Frobenius torus \(T_{x}(\rho)\),
 and such that \(\cp{F_{x,\rho}^{n}}\) has coefficients in~\(E\).
 Let \((\alpha_{i})_{i}\) be the roots of \(\cp{F_{x,\rho}^{n}}\).
 Let \(\Gamma \subset \bar{E}^{\star}\) be the subgroup
 generated by the \(\alpha_{i}\);
 it is a free abelian group that may be canonically identified with
 the character lattice \(\Hom(T_{x}(\rho),\Gm[E_{\lambda}])\).
 Let \(T\) be the algebraic torus over~\(E\)
 whose character lattice is~\(\Gamma\).
 By construction we have \(T_{E_{\lambda}} \cong T_{x}(\rho)\).

 The upshot of this computation is that we may
 view \(T_{x}(\rho)\) in a canonical way as an algebraic torus over~\(E\),
 if \(\rho\) is an \(E\)-rational Galois representation.
\end{remark}

\begin{proposition} \label{compatisom} 
 Let \(K\) be a finitely generated field.
 Let \(E\) be a number field; and
 let \(\lambda\) be a finite place of~\(E\).
 For \(i = 1,2\),
 let \(\rho_{i}\) be a \(\lambda\)-adic Galois representation of~\(K\).
 If \(\rho_{1}\) and~\(\rho_{2}\) are semisimple, quasi-compatible,
 and \(\GG_{\lambda}(\rho_{1} \oplus \rho_{2})\) is connected,
 then \(\rho_{1} \cong \rho_{2}\).
 \begin{proof}
  See \cref{densprop} through \cref{compatisomprf}.
 \end{proof}
\end{proposition}

\begin{theorem} \label{CSRisom} 
 Let \(K\) be a finitely generated field.
 Let \(E\) be a number field; and
 let \(\Lambda\) be a set of finite places of~\(E\).
 Let \(\rho_{\Lambda}\) and \(\rho'_{\Lambda}\)
 be two quasi-compatible systems of semisimple Galois representations.
 Assume that \(\GG_{\lambda}(\rho_{\lambda} \oplus \rho'_{\lambda})\)
 is connected for all \(\lambda \in \Lambda\).
 If there is a \(\lambda \in \Lambda\)
 such that \(\rho_{\lambda} \cong \rho'_{\lambda}\),
 then \(\rho_{\Lambda} \cong \rho'_{\Lambda}\).
 \begin{proof}
  This is an immediate consequence of \cref{compatisom}.
 \end{proof}
\end{theorem}

\begin{proposition} \label{recoverEnd} 
 Let \(K\) be a finitely generated field.
 Let \(E\) be a number field,
 and let \(\Lambda\) be the set of finite places of~\(E\)
 whose residue characteristic is different from~\(\chrc(K)\).
 Let \(\primes\) be the set of prime numbers different from~\(\chrc(K)\).
 Let \(\rho_{\Lambda}\) be a
 quasi-compatible system of semisimple Galois representations of~\(K\).
 Let \(\rho_{\primes}\) be the quasi-compatible system
 of Galois representations obtained by restricting to \(\QQ \subset E\),
 as in~\cref{restrictE};
 in other words, \(\rho_{\ell} = \bigoplus_{\lambda | \ell} \rho_{\lambda}\).
 Assume that \(\Gl(\rho_{\ell})\) is connected for all \(\ell \in \primes\).
 Fix \(\lambda_{0} \in \Lambda\).
 Define the field \(E' \subset E\) to be the
 subfield of~\(E\) generated by elements \(e \in E\)
 that satisfy the following condition:

 {\narrower\noindent
  There exists a model~\(X\) of~\(K\),
  a point \(x \in X^{\cl}\),
  and an integer~\(n \ge 1\),\\
  such that \(P_{x,\rho_{\lambda_{0}},n}(t) \in E[t]\)
  and \(e\) is a coefficient of \(P_{x,\rho_{\lambda_{0}},n}(t)\).
  \par}

 \noindent
 Let \(\ell\) be a prime number that splits completely in~\(E/\QQ\).
 If \(\End_{\Gal(\bar{K}/K),\QQl}(\rho_{\ell}) \cong E \otimes \QQl\),
 then \(E = E'\).
 \begin{proof}
  We restrict our attention to a finite subset of~\(\Lambda\),
  namely \(\Lambda_{0} = \{\lambda_{0}\} \cup \{ \lambda | \ell \}\).
  Let \(U \subset X\) be an open subset such that
  for all \(\lambda_{1}, \lambda_{2} \in \Lambda\)
  the representations
  \(\rho_{\lambda_{1}}\) and~\(\rho_{\lambda_{2}}\)
  are quasi-compatible at all \(x \in U^{\cl}\).
  For each \(x \in U^{\cl}\),
  let \(n_{x}\) be an integer such that
  \(P_{x}(t) = P_{x,\rho_{\lambda},n_{x}}(t) \in E[t]\)
  does not depend on \(\lambda \in \Lambda_{0}\).

  Let \(\lambda'\) be a place of \(E'\) above~\(\ell\).
  Let \(\lambda_{1}\) and~\(\lambda_{2}\) be two places
  of~\(E\) that lie above \(\lambda'\).
  We view \(\rho_{\lambda_{1}}\) and~\(\rho_{\lambda_{2}}\)
  as \(\lambda'\)-adic representations.
  Since \(\ell\) splits completely in~\(E/\QQ\),
  the inclusions \(\QQl \subset E'_{\lambda'} \subset E_{\lambda_{i}}\)
  are isomorphisms.
  By definition of~\(E'\) we have \(P_{x}(t) \in E'[t]\).
  Therefore \(\rho_{\lambda_{1}}\) and~\(\rho_{\lambda_{2}}\)
  are quasi-compatible \(\lambda'\)-adic representations;
  hence they are isomorphic by \cref{compatisom}.
  Let \(\rho_{\lambda'}\) be the \(\lambda'\)-adic Galois representation
  \(\bigoplus_{\lambda | \lambda'} \rho_{\lambda}\),
  as in~\cref{restrictE}.
  We conclude that
  \(\End_{\Gal(\bar{K}/K),E'_{\lambda'}}(\rho_{\lambda'})
   \cong \Mat_{[E:E']}(E'_{\lambda'})\),
  which implies \([E:E'] = 1\).
 \end{proof}
\end{proposition}

\paragraph{} \label{densprop} 
Let \(K\) be a finitely generated field,
and let \(X\) be a model of~\(K\).
There is a good notion of density for subsets of~\(X^{\cl}\).
This is described by Serre in~\cite{Se65} and~\cite{NXp},
and by Pink in appendix~B of~\cite{Pi97}.
For the convenience of the reader, we list some features of these densities.
Most of the following list is a reproduction of the statement of
proposition~B.7 of~\cite{Pi97}.
Let \(T \subset X^{\cl}\) be a subset.
If \(T\) has a density, we denote it with~\(\mu_{X}(T)\).
\begin{enumerate}
 \item If \(T \subset X^{\cl}\) has a density,
  then \(0 \le \mu_{X}(T) \le 1\).
 \item The set \(X^{\cl}\) has density~\(1\).
 \item If \(T\) is contained in a proper closed subset of~\(X\),
  then \(T\) has density~\(0\).
 \item If \(T_{1} \subset T \subset T_{2} \subset X^{\cl}\)
  such that \(\mu_{X}(T_{1})\) and \(\mu_{X}(T_{2})\) exist and are equal,
  then \(\mu_{X}(T)\) exists and is equal to
  \(\mu_{X}(T_{1}) = \mu_{X}(T_{2})\).
 \item \label{denscupcap} If \(T_{1}, T_{2} \subset X^{\cl}\) are two subsets,
  and three of the following densities exist,
  then so does the fourth, and we have
  \[
   \mu_{X}(T_{1} \cup T_{2}) + \mu_{X}(T_{1} \cap T_{2})
   = \mu_{X}(T_{1}) + \mu_{X}(T_{2}).
  \]
 \item If \(u \colon X \to X'\) is a birational morphism,
  then \(T\) has a density if and only if \(u(T)\) has a density,
  and if this is the case, then \(\mu_{X}(T) = \mu_{X'}(u(T))\).
\end{enumerate}

\paragraph{} 
Chebotarev's density theorem generalises to this setting.
Let \(Y \to X\) be a finite \'{e}tale Galois covering
of integral schemes of finite type over \(\Spec(\ZZ)\).
Denote the Galois group with~\(G\).
For each point \(y \in Y^{\cl}\) with image \(x \in X^{\cl}\)
the inverse of the Frobenius endomorphism of \(\kappa(y)/\kappa(x)\)
determines an element \(F_{y} \in G\).
The conjugacy class of \(F_{y}\) only depends on~\(x\),
and we denote it with~\(\mathcal{F}_{x}\).

\begin{theorem} \label{Cheb} 
 Let \(Y \to X\) be a finite \'{e}tale Galois covering
 of integral schemes of finite type over \(\Spec(\ZZ)\) with group~\(G\).
 For every conjugacy class \(C \subset G\),
 the set \(\{ x \in X^{\cl} \mid \mathcal{F}_{x} = C \}\)
 has density \(\frac{\#C}{\#G}\).
 \begin{proof}
  See proposition~B.9 of~\cite{Pi97}.
 \end{proof}
\end{theorem}

\begin{theorem} 
 Let \(K\) be a finitely generated field.
 Let \(E\) be a number field, and let \(\lambda\) be a finite place of~\(E\).
 Let \(\rho\) be a semisimple \(\lambda\)-adic Galois representation of~\(K\).
 Assume that \(\GG_{\lambda}(\rho)\) is connected.
 There is a non-empty Zariski open subset \(U \subset \GG_{\lambda}(\rho)\)
 such that for every model~\(X\) of~\(K\),
 and every closed point \(x \in X^{\cl}\),
 if \(\rho\) is unramified at~\(x\),
 and for some \(n \ge 1\)
 the Frobenius element~\(F_{x,\rho}^{n}\)
 is conjugate to an element of~\(U(E_{\lambda})\),
 then \(T_{x}(\rho)\) is a maximal torus of~\(\GG_{\lambda}(\rho)\).
 \begin{proof}
  See theorem~3.7 of~\cite{Ch92}.
  The statement in~\cite{Ch92} is for abelian varieties,
  but the proof is completely general.
 \end{proof}
\end{theorem}

\begin{corollary}[3.8 of~\cite{Ch92}] \label{maxtor} 
 Let \(K\) be a finitely generated field.
 Let \(E\) be a number field, and let \(\lambda\) be a finite place of~\(E\).
 Let \(\rho\) be a semisimple \(\lambda\)-adic Galois representation of~\(K\).
 Assume that \(\GG_{\lambda}(\rho)\) is connected.
 Let \(X\) be a model of~\(K\).
 Let \(\Sigma \subset X^{\cl}\) be the set of points \(x \in X^{\cl}\)
 for which \(\rho\) is unramified at~\(x\)
 and \(T_{x}(\rho)\) is a maximal torus of \(\GG_{\lambda}(\rho)\).
 Then \(\Sigma\) has density~\(1\).
\end{corollary}

\begin{lemma} \label{BNalggrp} 
 Let \(K\) be a field of characteristic~\(0\).
 Let \(G\) be a reductive group over~\(K\).
 Let \(\rho_{1}\) and~\(\rho_{2}\) be two
 finite-dimensional semisimple representations of~\(G\).
 Let \(S \subset G(K)\) be a subset that is Zariski-dense in~\(G\).
 Assume that for all \(g \in S\) we have
 \(\trace(\rho_{1}(g)) = \trace(\rho_{2}(g))\).
 Then \(\rho_{1} \cong \rho_{2}\) as representations of~\(G\).
 \begin{proof}
  Note that \(\trace \circ \rho_{i}\) is
  a separated morphism of schemes \(G \to \mathbb{A}^{1}_{K}\).
  Therefore we have \(\trace(\rho_{1}(g)) = \trace(\rho_{2}(g))\)
  for all \(g \in G(K)\).
  By linearity, we find that
  \(\trace(\rho_{1}(\alpha)) = \trace(\rho_{2}(\alpha))\)
  for all \(\alpha\) in the group algebra \(K[G(K)]\).
  By proposition~3{} in \S12, \textnumero1 of~\cite{BouAlg},
  we conclude that \(\rho_{1} \cong \rho_{2}\) as representations of \(G(K)\),
  hence as representations of~\(G\).
 \end{proof}
\end{lemma}

\begin{lemma} \label{maxtorisom} 
 Let \(K\) be a finitely generated field.
 Let \(E\) be a number field; and
 let \(\lambda\) be a finite place of~\(E\).
 For \(i = 1,2\),
 let \(\rho_{i}\) be
 a semisimple \(\lambda\)-adic Galois representation of~\(K\).
 Write \(\rho\) for \(\rho_{1} \oplus \rho_{2}\).
 Assume that \(\GG_{\lambda}(\rho)\) is connected.
 If there is a model~\(X\) of~\(K\),
 and a point \(x \in X^{\cl}\)
 such that
 \(\rho\) is unramified at~\(x\),
 and \(T_{x}(\rho)\) is a maximal torus,
 and \(P_{x,\rho_{1},n}(t) = P_{x,\rho_{2},n}(t)\)
 for some \(n \ge 1\),
 then \(\rho_{1} \cong \rho_{2}\) as \(\lambda\)-adic Galois representations.
 \begin{proof}
  Write \(T\) for \(T_{x}(\rho)\).
  Observe that
  \(P_{x,\rho_{1},kn}(t) = P_{x,\rho_{2},kn}(t)\)
  for all \(k \ge 1\).
  Let \(H_{n}\) be the algebraic subgroup of~\(\GG_{\lambda}(\rho)\)
  that is generated by~\(F_{x,\rho}^{n}\).
  Recall that \(T\) is the identity component of~\(H_{n}\).
  Note that for some \(k \ge 1\),
  we have \(F_{x,\rho}^{kn} \in T(E_{\lambda})\).
  Replace \(n\) by \(kn\), so that we may assume
  that \(F_{x,\rho}^{n}\) generates~\(T\) as algebraic group.

  The set \(\{ F_{x,\rho}^{kn} \mid k \ge 1 \}\)
  is a Zariski dense subset of~\(T\).
  Since \(P_{x,\rho_{1},kn}(t) = P_{x,\rho_{2},kn}(t)\),
  for all \(k \ge 1\),
  \cref{BNalggrp} implies that
  \(\rho_{1}|_{T} \cong \rho_{2}|_{T}\).
  Because \(T\) is a maximal torus of~\(\GG_{\lambda}(\rho)\)
  and \(\GG_{\lambda}(\rho)\) is connected,
  we find that \(\rho_{1} \cong \rho_{2}\) as representations
  of~\(\GG_{\lambda}(\rho)\),
  and hence as \(\lambda\)-adic Galois representations of~\(K\).
 \end{proof}
\end{lemma}

\begin{nproof}[of \cref{compatisom}] \label{compatisomprf} 
 Let \(X\) be a model of~\(K\).
 By \cref{maxtor}, the subset of points \(x \in X^{\cl}\)
 for which \(T_{x}(\rho_{1} \oplus \rho_{2})\)
 is a maximal torus is a subset with density~\(1\).
 By definition of compatibility, and \cref{densprop},
 the subset of points \(x \in X^{\cl}\)
 at which \(\rho_{1}\) and \(\rho_{2}\) are quasi-compatible
 is also a subset with density~\(1\).
 By \cref{denscupcap} these subsets have non-empty intersection:
 there exists a point \(x \in X^{\cl}\)
 such that \(T_{x}(\rho_{1} \oplus \rho_{2})\) is a maximal torus
 and \(\rho_{1}\) and \(\rho_{2}\) are quasi-compatible at~\(x\).
 Now \cref{compatisom} follows from \cref{maxtorisom}.
\end{nproof}

\begin{lemma} \label{rkindeplambda} 
 Let \(K\) be a finitely generated field.
 Let \(E\) be a number field,
 and let \(\Lambda\) be a set of finite places of~\(E\).
 Let \(\rho_{\Lambda}\) be
 a quasi-compatible system of semisimple Galois representation of~\(K\).
 Assume that \(\GG_{\lambda}(\rho_{\lambda})\) is connected,
 for all \(\lambda \in \Lambda\).
 Then the absolute rank of \(\GG_{\lambda}(\rho_{\lambda})\)
 is independent of~\(\lambda\).
 \begin{proof}
  It suffices to assume that \(\Lambda = \{\lambda_{1}, \lambda_{2}\}\).
  Let \(X\) be a model of~\(K\).
  For \(i = 1,2\),
  let \(\Sigma_{i} \subset X^{\cl}\) be the set of points \(x \in X^{\cl}\)
  for which \(\rho_{\lambda_{i}}\) is unramified at~\(x\)
  and \(T_{x}(\rho_{\lambda_{i}})\) is a maximal torus
  of \(\GG_{\lambda_{i}}(\rho_{\lambda_{i}})\).
  Then \(\Sigma_{i}\) has density~\(1\), by \cref{maxtor}.
  Let \(U \subset X\) be an open subset
  such that \(\rho_{\lambda_{1}}\) and~\(\rho_{\lambda_{2}}\)
  are compatible at all \(x \in U^{\cl}\).

  Put \(\Sigma = \Sigma_{1} \cap \Sigma_{2} \cap U^{\cl}\);
  by \cref{denscupcap} we know that \(\Sigma\) is non-empty.
  Fix a closed point \(x \in \Sigma\).
  Since \(\rho_{\lambda_{1}}\) and~\(\rho_{\lambda_{2}}\)
  are \(E\)-rational and quasi-compatible at~\(x\),
  there exists a torus~\(T\) over~\(E\)
  such that \(T_{E_{\lambda_{i}}} \cong T_{x}(\rho_{\lambda_{i}})\),
  see \cref{FrobTorE}.
  The tori \(T_{x}(\rho_{\lambda_{i}})
   \subset \GG_{\lambda_{i}}(\rho_{\lambda_{i}})\)
  are maximal tori, by assumption.
  We conclude that \(\GG_{\lambda_{1}}(\rho_{\lambda_{1}})\)
  and~\(\GG_{\lambda_{2}}(\rho_{\lambda_{2}})\) have the same absolute rank.
 \end{proof}
\end{lemma}

\section{Remark on the Mumford--Tate conjecture} \label{MTCrmk} 

\readme 
In this section we recall the Mumford--Tate conjecture.
A priori, this conjecture depends on the choice of a prime number~\(\ell\).
We show that for abelian motives this conjecture does not depend on~\(\ell\)
(see \cref{MTCindepl}).

\begin{conjecture}[Mumford--Tate] \label{MTC} 
 Let \(M\) be a motive over a finitely generated subfield of~\(\CC\).
 Let \(\ell\) be a prime number.
 Under the comparison isomorphism
 \(\HB(M) \otimes \QQl \cong \Hl(M)\),
 see \cref{Artin}, we have
 \[
  \GB(M) \otimes \QQl \cong \Glc(M).
 \]
\end{conjecture}

\paragraph{} 
If \(M\) is a motive over a finitely generated subfield of~\(\CC\), then
we write \(\MTC(M)_{\ell}\) for the conjectural statement
\(
 \GB(M) \otimes \QQl \cong \Glc(M).
\)

\begin{lemma} \label{MTCext} 
 Let \(K \subset L\) be finitely generated subfields of~\(\CC\).
 Let \(M\) be a motive over~\(K\).
 Then \(\GB(M_{L}) = \GB(M)\) and \(\Glc(M_{L}) = \Glc(M)\).
 In particular \(\MTC_{\ell}(M) \iff \MTC_{\ell}(M_{L})\).
 \begin{proof}
  See proposition~1.3 of~\cite{Mo15}.
 \end{proof}
\end{lemma}

\begin{proposition} \label{MTCcentres} 
	Let \(M\) be an abelian motive
 over a finitely generated field \(K \subset \CC\).
	Let \(\ZB(M)\) be the centre of the Mumford--Tate group~\(\GB(M)\),
 and let \(\Zl(M)\) be the centre of~\(\Glc(M)\).
	Then \(\Zl(M) \subset \ZB(M) \otimes \QQl\),
 and \(\Zl(M)^{\circ} = \ZB(M)^{\circ} \otimes \QQl\).
	\begin{proof}
		The result is true for abelian varieties
		(see theorem~1.3.1 of~\cite{Va08} or corollary~2.11 of~\cite{UY13}).

		By definition of abelian motive,
  there is an abelian variety~\(A\) such that \(M\) is
  contained in the Tannakian subcategory of motives generated by~\(\HH(A)\).
  By \cref{MTCext}, we may assume that \(\Gl(A)\) is connected
  for all prime numbers~\(\ell\), without loss of generality.

		We have a surjection \(\GB(A) \onto \GB(M)\).
  Since \(\GB(A)\) is reductive,
  \(\ZB(M)\) is the image of \(\ZB(A)\) under this map.
		The same is true on the \(\ell\)-adic side.
  (Note that \(\Gl(A)\) is reductive, by Satz~3{} in~\S5 of~\cite{Fal83};
  see also~\cite{Fal84}.)
		Thus we obtain a commutative diagram with solid arrows
		\[
			\begin{tikzcd}
    \Zl(A) \dar[hook] \rar[two heads] &
    \Zl(M) \dar[dotted,hook] \rar[hook] &
				\Gl(M) \dar[hook] \\
    \ZB(A) \otimes \QQl \rar[two heads] &
				\ZB(M) \otimes \QQl \rar[hook] &
				\GB(M) \otimes \QQl
			\end{tikzcd}
		\]
		which shows that the dotted arrow exists and is an inclusion.
  (The vertical arrow on the right exists and is an inclusion,
  by \cref{hdgmot}.)

  Finally, observe that \(\ZB(M)\) and~\(\Zl(M)\) have the same rank.
  Indeed, as remarked at the beginning of the proof,
  we know that \(\Zl(A)^{\circ} \into \ZB(A)^{\circ} \otimes \QQl\)
  is an isomorphism.
  The commutative diagram above shows that the inclusion
  \(\Zl(M)^{\circ} \into \ZB(M)^{\circ} \otimes \QQl\)
  must be an isomorphism.
	\end{proof}
\end{proposition}

\begin{proposition} \label{rkeqimplMTC} 
 Let \(M\) be an abelian motive
 over a finitely generated subfield \(K \subset \CC\).
 Let \(\ell\) be a prime number.
 If \(\GB(M)\) and~\(\Gl(M)\) have the same absolute rank,
 then \(\MTC(M)_{\ell}\) is true.
 \begin{proof}
  We apply the Borel--de Siebenthal theorem
  (see~\cite{BdS}; or~\cite{PL2015}):
  since \(\Glc(M) \subset \GB(M) \otimes \QQl\) has maximal rank,
  it is equal to the connected component of the centraliser of its centre.
  By \cref{MTCcentres}, we know that the centre of \(\Glc(M)\)
  is contained in the centre of \(\GB(M) \otimes \QQl\).
  Hence \(\Glc(M) \cong \GB(M) \otimes \QQl\).
 \end{proof}
\end{proposition}

\begin{corollary} \label{MTCindepl} 
 Let \(M\) be an abelian motive
 over a finitely generated subfield \(K \subset \CC\).
 The Mumford--Tate conjecture is independent of the
 choice of the prime number~\(\ell\).
 \begin{proof}
  Without loss of generality,
  we may and do assume that the groups \(\Gl(M)\)
  are connected for all prime numbers~\(\ell\), by \cref{MTCext}.
  If the Mumford--Tate conjecture is true for one prime number~\(\ell\),
  then the groups \(\GB(M)\) and \(\Gl(M)\) have the same absolute rank.
  By \cref{abmotCSR} the Galois representations~\(\Hl(M)\)
  form a quasi-compatible system,
  and by \cref{rkindeplambda} the rank of the groups \(\Gl(M)\)
  does not depend on~\(\ell\).
  The result follows from \cref{rkeqimplMTC}
 \end{proof}
\end{corollary}

\printbibliography
\end{document}